\documentclass[11pt,A4paper]{article}

\usepackage[pdftex,hmargin={2cm,2cm},vmargin={2cm,2cm}]{geometry}
\usepackage[utf8]{inputenc}
\usepackage{amsmath}
\usepackage{amsfonts}
\usepackage{graphicx} 
\usepackage[table]{xcolor}
\usepackage[pdftex,bookmarks=false,colorlinks=true,citecolor=blue]{hyperref}

\newcommand{\R}{\mathbb{R}}
\newcommand{\RR}{\mathcal{R}}
\newcommand{\N}{\mathbb{N}}
\newcommand{\C}{\mathbb{C}}
\newcommand{\G}{\mathcal{G}}
\renewcommand{\d}{{\rm d}}
\newcommand{\ep}{\varepsilon}
\newcommand{\Co}{\mathcal{C}}
\renewcommand{\L}{\mathcal{L}}
\newcommand{\bif}{\textnormal{bif}}
\newcommand{\X}{\mathcal{X}}

\newcommand{\Tau}{\mathcal{T}}
\newcommand{\diag}{\textnormal{diag}}
\newcommand{\Zz}{\mathcal{Z}}
\newcommand{\ess}{\textnormal{ess}}
\newcommand{\TT}{\mathbf{T}}

\newtheorem{theorem}{Theorem}[section]

\newtheorem{assumption}[theorem]{Assumption}

\newtheorem{lemma}[theorem]{Lemma}

\date{}
\author{Quentin Richard$^a$, Marc Choisy$^{b,c}$, Thierry Lefèvre$^a$ and Ramsès Djidjou-Demasse$^{a,*}$}

\date{	{\small $^{a}$MIVEGEC, Univ. Montpellier, IRD, CNRS, Montpellier, France} \\
	{\small $^{b}$ Centre for Tropical Medicine and Global Health, Nuffield Department of Medicine, University of Oxford, UK} \\
	{\small $^{c}$ Oxford University Clinical Research Unit, Ho Chi Minh City, Vietnam} \\
	{\small $^*$Author for correspondence: ramses.djidjoudemasse@ird.fr}}

\title{Human-vector malaria transmission model structured by age, time since infection and waning immunity}

\begin{document}

\maketitle

\begin{abstract}
In contrast to the many theoretical studies on the transmission of human-mosquitoes malaria infection, few studies have considered a multiple structure model formulations including (i) the chronological age of humans and mosquitoes population, (ii) the time since humans and mosquitoes  are infected and (iii) humans waning immunity ({\it i.e.,} the progressive loss of protective antibodies after recovery). Such structural variables are well documented to be fundamental for the transmission of human-mosquitoes malaria infections. Here we formulate an age-structured model accounting for the three structural variables. Using integrated semigroups theory, we first handle the well-posedness of the model proposed. We also investigate the existence of model's steady-states. A disease-free equilibrium always exists while the existence of endemic equilibria is discussed. We derive the threshold $\RR_0$ (the basic reproduction number). The expression of the $\RR_0$ obtained here particularly highlight the effect of above structural variables on key important epidemiological traits of the human-vector association. This includes, humans and mosquitoes transmission probability and survival rates. Next, we derive a necessary and sufficient condition that implies the bifurcation of an endemic equilibrium. In some configuration where the age-structure of the human population is neglected,  we show that, depending on the sign of some constant $C_{\bif}$ given by the parameters, a bifurcation occurs at $\RR_0=1$ that is either forward or backward. In the former case, it means that there exists a (unique) endemic equilibrium if and only if $\RR_0>1$. In the latter case, no endemic equilibrium exists for $\RR_0\ll 1$ small enough, a unique exists if $\RR_0>1$ while multiple endemic equilibria exist when $0\ll \RR_0<1$ close enough to $1$.

\vspace{0.2in}\noindent \textbf{Key words}. Vector-borne diseases; Malaria; Basic reproduction number; Age-structured model; Bifurcation analysis.
\end{abstract}

\section{Introduction}

More than one century ago, Ross \cite{Ross1911} introduced the first mathematical model for the transmission of malaria. The latter model was refined later by MacDonald \cite{MacDonald1957}. This vector-borne disease is still a wide subject of study in epidemiology, see \textit{e.g.} \cite{Brauer2019} and the references therein. Most mathematical models about the transmission dynamic of malaria use ordinary differential equations \cite{CaiMartcheva2017,ChitnisCushing2006,ChitnisCushing2008,Djidjou2020,DucrotZongo2009,GaoRuan2014,Harvim2019,LiMartcheva2015,LouZhao2010,PropserMartcheva2014,TchoumiSallet2018,Traore2020,Traore2017,Xing2020}.

Besides the age dependence of the natural mortality rate of the human population, age also plays an important role in the transmission dynamic of malaria. Indeed, while more than 400 000 deaths occurred in 2019 due to malaria infections, about 67$\%$ were among the youngest population, {\it i.e.} less than 5 years old \cite{WHO2020}. Furthemore, it is becoming clear that the human infectious reservoir is also age-dependant, with 5 to 15 year-old children representing the most important source of infection to mosquitoes \cite{Coalson2018,Felger2012}. It is then crucial to take into account an age-structure in the host population as in \cite{Aguas2008,Geisse2012,MaWarner2019}.

Another key factor is the time since a human host is infected. It is particularly relevant since  the production of gametocytes (sexual forms of malaria parasites) within a human host is strongly related to the time since the host is infected \cite{djidjoudemasse2020understanding}. Moreover, there is a clear relationship between gametocyte density and the transmission probability per bite from human to mosquito \cite{ChurcherEtAl2013}. 

Recent works \cite{Bellan2010,Styer2007} emphasized the importance of mosquitoes senescence as part of the modelling procedure in mosquito populations (see also \cite{LefevreWork}). In the literature, considering chronological age-structured mosquitoes was only taken into account in \cite{Rock2015}. This consideration is important since mosquitoes live on average 14 days \cite{ChitnisCushing2008,Rock2015}, while the extrinsic incubation period \cite{Ohm2018} is in average 11-12 days \cite{ChitnisCushing2008}. Hence a mosquito that become infected at the end of its lifespan, will probably never infect any human. The probability of transmission from mosquitoes to humans should consequently depend both on the chronological and infection ages of mosquitoes.

In the literature, age-structured models may incorporate additional structures such as body size, space or more general phenotypic trait, see \textit{e.g.} \cite{Djidjou2017,KangRuan2020,Sinko67,Thieme91,Webb2008} and the references therein for a survey of such models. However, epidemiological models including both infection and chronological age are not so common in the literature (see \cite{BurieDucrot2017,Dietz85,Hoppensteadt74,Inaba2006b,Inaba2016,Kapitanov2015,LarochePerasso2016,Richard2020,Zhou2002} for an exhaustive list). In a context of vector-borne infectious diseases, and more precisely focusing on malaria transmission, models with both chronological and infection age structures have never been considered until now. The present model takes into account these structures simultaneously in humans and mosquitoes populations. Moreover, we also consider the time since an human recovered as another continuous variable to account for a potential waning immunity ({\it i.e.,} the progressive loss of protective antibodies after recovery).

In this paper, we first handle the well-posedness of the model. To this end, we use integrated semigroups theory, whose approach was introduced in \cite{Arendt87,PratoSinestrari87,KellermanHieber89,Neubrander88,Thieme90}. We also refer to \cite{MagalRuan2018} and the references therein for more details. Note that this framework was successfully used \cite{BurieDucrot2017} in a context of a population with double structure. However, in our case, the shape of the force of infection and more precisely the fact that it has a singularity when the total population of humans is zero, makes the analysis much more delicate. We use a classical fixed point argument in an appropriate $L^1$ space combined with some estimates of the populations. One could also proceed with the classical method, that is, use solutions integrated along the characteristics and work with nonlinear Volterra equations. We refer to the monographs \cite{Iannelli94,DiekmannMetz86,Webb85} on this method. 

We also investigate the existence of steady-states, that are time-independent solutions of the model. A disease-free equilibrium clearly always exists while the existence of endemic equilibria is discussed. We derive the threshold $\RR_0$ (the basic reproduction number) and a necessary and sufficient condition that implies the existence of an endemic equilibrium. While it is difficult to exploit in the general case, we then focus on a particular case where the age-structure of the human population is neglected, and the latter condition becomes explicit. We show that, depending on the sign of a constant $C_{\bif}$ given by the parameters, a bifurcation occurs at $\RR_0=1$ that is either forward or backward (see \textit{e.g.} \cite{ChitnisCushing2006,DucrotZongo2009,Inaba2017,TchoumiSallet2018} for more details on such bifurcations). In the former case, it means that there exists a (unique) endemic equilibrium if and only if $\RR_0>1$. In the latter case, no endemic equilibrium exists for $\RR_0\ll 1$ small enough, a unique exists if $\RR_0>1$ while multiple endemic equilibria exist when $0\ll \RR_0<1$ close enough to $1$.

The integrated semigroups framework allows us to linearize the system around each equilibrium and obtain linear $\Co_0$-semigroups. Using spectral theory, we are able to prove the local stability of the disease-free equilibrium under the condition $\RR_0<1$, while it is unstable whenever $\RR_0>1$ (see \textit{e.g.} \cite{EngelNagel2000,Webb85} for more results on this topic).

The paper is structured as follows: we first introduce the model and define the outputs and models parameters. Next, we state and discuss main results that will be obtained in this work. These include the existence and uniqueness of bounded solutions, the disease invasion process and the bifurcation --forward and backward-- of an endemic equilibrium in some special cases. Numerical simulations are provided to illustrate above main results. Finally, details on the proof of the main results of this work complete the paper.


\section{Description of the model} \label{Sec:model}

\subsection{Model overview}

At time $t\geq 0$, the density of humans with age $a\geq 0$, that are susceptible to the infection is denoted by $S_h(t,a)$. These individuals can become infected due to bites of infected mosquitoes with a rate $\lambda_{m\to h}(t,a)$, called the force of infection of mosquitoes to humans with age $a$. Infected humans population is additionally structured by the time since infection, so that $I_h(t,a,\tau)$ denotes the density at time $t$ of individuals of age $a$ that have been infected for a duration $\tau \geq 0$. During their infection, humans can either recover at a time since infection $\tau$ with rate $\gamma_h(a,\tau)$, or die from the infection with the rate $\nu_h(a,\tau)$. At time $t$, humans $R(t,a,\eta)$ with age $a$, that have recovered from the infection for a duration $\eta$, are temporarily immunized and lose their immunity at rate $k_h(a,\eta)$. Each human may also die due to natural causes with an age-dependent rate $\mu_h(a)$. The flux of newborn humans is assumed constant to $\Lambda_h$. 

\begin{figure}[!h]
\begin{center}\includegraphics[width=.7\linewidth]{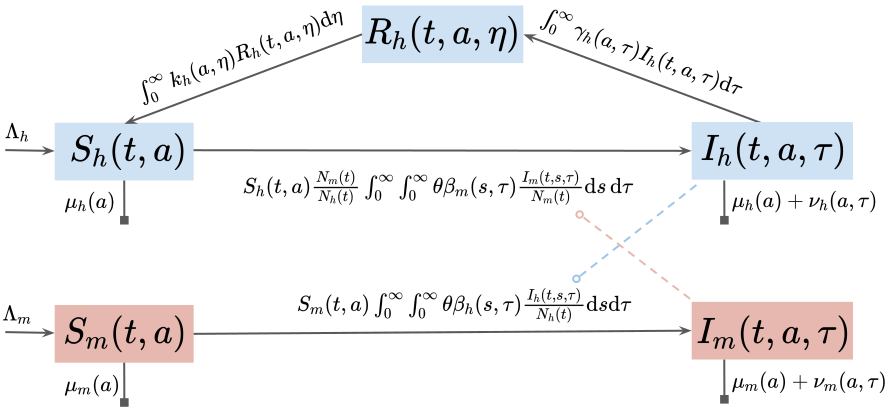}\end{center}
\caption{{\bf Flow diagram illustrating the interactions between humans (subscript $h$) and mosquitoes (subscript $m$).} Newborns humans are recruited at a constant rate $\Lambda_h$. Natural death rate of humans aged $a$ is $\mu_h(a)$, and if infected since time $\tau$, the disease induced mortality is $\nu_h(a,\tau)$. The force of infection from mosquitoes to humans at time $t$ --- $\frac{N_m(t)}{N_h(t)}\int_0^\infty \int_0^\infty \theta \beta_m(s,\tau)\frac{I_m(t,s,\tau)}{N_m(t)}\d s~\d\tau$ --- is defined as the product of the probability that mosquitoes aged $s$ and infected since time $\tau$ remain infectious $I_m(t,s,\tau)/N_m(t)$, the number of mosquito bites per human per time $\theta N_m(t)/N_h(t)$ and the probability of disease transmission from mosquitoes to humans $\beta_m(s,\tau)$. Humans aged $s$ and infected since time $\tau$ recover from the disease at rate $\gamma_h(s,\tau)$, get temporary infection immunity which they lose at rate $k_h(s,\eta)$ after a duration $\eta$ of recovery. Newborns mosquitoes are generated at rate $\Lambda_m$. Natural death rate of mosquitoes aged $a$ is $\mu_m(a)$, and if infected since time $\tau$, the disease induced mortality is $\nu_m(a,\tau)$. The force on infection from humans to mosquitoes at time $t$ --- $ \int_0^\infty \int_0^\infty \theta \beta_h(s,\tau)\frac{I_h(t,s,\tau)}{N_h(t)}\d s~\d\tau$ --- is defined as the product of the probability of disease transmission from human aged $s$ and infected since time $\tau$ to the mosquito
$\beta_h(s,\tau)$, mosquitoes biting rate $\theta$ and the
probability that  human is infectious $I_h(t,s,\tau)/N_h(t)$.} \label{Fig:diagram}
\end{figure}

At time $t$, susceptible mosquitoes of age $a$, denoted by $S_m(t,a)$, become infected by taking contaminated blood from infected humans at rate $\lambda_{h\to m}(t,a)$, called the force of infection of humans to mosquitoes with age $a$. The mosquito population that have been infected for a duration $\tau$, $I_m(t,a,\tau)$, may die from the infection at rate $\nu_m(a,\tau)$. As for the human population, the age-dependent natural death of the mosquitoes is $\mu_m(a)$ while the flux of newborn mosquitoes is assumed constant to $\Lambda_m$. The human-mosquitoes infection life cycle is shown is Figure \ref{Fig:diagram}. The total number of humans and mosquitoes at time $t$ are respectively given by
$$N_h(t)=\int_0^\infty S_h(t,a)\d a+\int_0^\infty \int_0^\infty I_h(t,a,\tau)\d a~\d\tau+\int_0^\infty \int_0^\infty R_h(t,a,\eta)\d a~\d \eta,$$
and
$$N_m(t)=\int_0^\infty S_m(t,a) \d a+\int_0^\infty \int_0^\infty I_m(t,a,\tau)\d a~\d\tau.$$

The force of infection from mosquitoes to humans with age $a$ is given by:
\begin{equation}\label{Eq:lambda_m}
\lambda_{m\to h}(t,a)=\frac{S_h(t,a)}{N_h(t)}\int_0^\infty \int_0^\infty \theta \beta_m(s,\tau)I_m(t,s,\tau)\d s~\d\tau.
\end{equation}
and describes the number of newly infected human with age $a$ at time $t$. It consists of the probability that human with age $a$ encountered by a mosquito is susceptible $S_h(t,a)/N_h(t)$ and the infection efficiency of the mosquito population $\int_0^\infty \int_0^\infty \theta \beta_m(s,\tau)I_m(t,s,\tau)\d s~\d\tau$. The latter efficiency takes into account (i) $\theta$ the number of human bitten by mosquito by unit of time and (ii) $\beta_m(s,\tau)$ the probability of disease transmission from one bite of infected mosquito (with age $a$ and which is infected since a time $\tau$) to a human. Similarly, the force of infection from human to mosquitoes with age $a$ is given by:
\begin{equation}\label{Eq:lambda_h}
    \lambda_{h\to m}(t,a)=\frac{S_m(t,a)}{N_h(t)} \int_0^\infty \int_0^\infty \theta \beta_h(s,\tau)I_h(t,s,\tau)\d s~\d\tau,
\end{equation}
where $\beta_h(s,\tau)$ is the probability of disease transmission from human with age $s$ (and infected since a time $\tau$) to a mosquito for each bite.

\subsection{The mathematical model}

Based on the above notations, the model considered reads as:
\begin{equation}\label{Eq:Model}
\left\{
\begin{array}{rcl}
\left(\frac{\partial}{\partial t}+\frac{\partial}{\partial a}\right)S_h(t,a)&=&\int_0^\infty k_h(a,\eta) R_h(t,a,\eta)\d \eta-\mu_h(a) S_h(t,a)-\lambda_{m\to h}(t,a) \\
\left(\frac{\partial}{\partial t}+\frac{\partial}{\partial a}+\frac{\partial}{\partial \tau}\right)I_h(t,a,\tau)&=&-\left(\mu_h(a)+\nu_h(a,\tau)+\gamma_h(a,\tau)\right)I_h(t,a,\tau), \vspace{0.1cm} \\
\left(\frac{\partial}{\partial t}+\frac{\partial}{\partial a}+
\frac{\partial}{\partial \eta}\right)R_h(t,a,\eta)&=&-(\mu_h(a)+k_h(a,\eta))R_h(t,a,\eta), \vspace{0.1cm} \\
\left(\frac{\partial}{\partial t}+\frac{\partial}{\partial a}\right)S_m(t,a)&=&-\mu_m(a)S_m(t,a)-\lambda_{h\to m}(t,a), \\
\left(\frac{\partial }{\partial t}+\frac{\partial }{\partial a}+\frac{\partial }{\partial \tau}\right)I_m(t,a,\tau)&=&-(\mu_m(a)+\nu_m(a,\tau))I_m(t,a,\tau), \\
\end{array}
\right.
\end{equation}
for each $(t,a,\tau,\eta)\in (0,\infty)^4$. System \eqref{Eq:Model} is associated to the following boundary conditions:
\begin{equation}\label{Eq:Bound_cond}
\left\{
\begin{array}{rclll}
S_h(t,0)&=&\Lambda_h, &S_m(t,0)=\Lambda_m,\\
I_h(t,a,0)&=&\lambda_{m\to h}(t,a), &I_h(t,0,\tau)=0, \vspace{0.1cm} \\
R_h(t,a,0)&=&\int_0^\infty \gamma_h(a,\tau)I_h(t,a,\tau) \d \tau, &R_h(t,0,\eta)=0, \vspace{0.1cm} \\
I_m(t,a,0)&=&\lambda_{h\to m}(t,a), &I_m(t,0,\tau)=0
\end{array}
\right.
\end{equation}
and the initial conditions (at $t=0$):
\begin{equation}\label{Eq:Ini_cond}
\left\{
\begin{array}{rclrlrl}
S_h(0,a)&=&S_{h,0}(a), &&I_h(0,a,\tau)&=&I_{h,0}(a,\tau), \qquad R_h(0,a,\eta)=R_{h,0}(a,\eta), \\
S_m(0,a)&=&S_{m,0}(a), &&I_m(0,a,\tau)&=& I_{m,0}(a,\tau), 
\end{array}
\right.
\end{equation}
for each $(a,\eta,\tau)\in \R_+^3$. The summary of all notations used in the latter model can be found in Table \ref{Tab-parameters}.

\begin{table}[!htp]
\begin{center}	\begin{small}
		\begin{tabular}{|c|l|l|}
		    \hline
			\textbf{Category} & \textbf{Description} & \textbf{Unit}\\
			\hline \hline
			\textbf{Notations} & & \\
			\hline
				$t$ & Time & Tu  \\
					$a$ & Chronological age & Tu \\
			$\tau$ & Time since infection  & Tu \\
					$\eta$ & Time since recovery for humans  & Tu \\
			\hline
			\textbf{States variables} & & \\
			\hline
			$S_h(t,a), S_m(t,a)$ & Susceptible humans and mosquitoes & No unit  \\ 
			$I_h(t,a,\tau), I_m(t,a,\tau)$ & Infected humans and mosquitoes & No unit \\
			$R_h(t,a,\eta)$ & Recovered humans & No unit \\ 
			$N_h(t), N_m(t)$ & Total humans and mosquitoes populations & No unit \\
			\hline 
    		\textbf{Parameters} & & \\
			\hline
			$\Lambda_h$ & Humans recruitment rate & Tu$^{-1}$ \\
		    $\Lambda_m$ & Mosquitoes recruitment rate & Tu$^{-1}$\\
        	$\mu_h(a)$ & Humans death rate & Tu$^{-1}$ \\
		    $\mu_m(a)$ & Mosquitoes death rate & Tu$^{-1}$ \\
        	$\nu_h(a,\tau)$ & Humans death rate induced by the infection & Tu$^{-1}$ \\
			$\nu_m(a,\tau)$ & Mosquitoes death rate induced by the infection & Tu$^{-1}$ \\
			$\gamma_h(a,\tau)$ & Recovery rate of humans infections & Tu$^{-1}$ \\
			$k_h(a,\eta)$& Rate at which humans lose their immunity & Tu$^{-1}$ \\
			$\beta_h(a,\tau)$ & Disease transmission probability from humans to mosquitoes & No unit \\
			$\beta_m(a,\tau)$ & Disease transmission probability from mosquitoes to humans & No unit \\
			$\theta$ & Biting rate of mosquitoes & Tu$^{-1}$ \\
				\hline \hline
		\end{tabular}\\
		Tu=time unit; h=human; m=mosquitoes
	\end{small}
	\end{center}
	\caption{Main notations, state variables and parameters of the model.} 
	\label{Tab-parameters}
\end{table}

\section{Main results}
In this section we will state the main results of this work. This include the existence and uniqueness of bounded solutions, the disease invasion process and the bifurcation of an endemic equilibrium in some special cases.

\subsection{Existence and uniqueness of bounded solutions}
Here we handle the well-posedness of \eqref{Eq:Model} with an integrated semigroups approach under the following general assumption 
\begin{assumption}\label{Assump}
\begin{enumerate}
    \item Recruitment rates $\Lambda_h$, $\Lambda_m$ and biting rate $\theta$ are positive constants;
    \item Mortality rates satisfy $\mu_h\in L^\infty_+(\R_+)$, $\nu_h\in L^\infty_+(\R_+^2)$, $\mu_m\in L^\infty_+(\R_+), \nu_m\in L^\infty_+(\R_+^2)$ and are such that there exists a constant $\mu_0>0$ so that $\mu_h(a)\geq \mu_0$ and $\mu_m(a)\geq \mu_0$ for each $a\in \R_+$;
    \item Transmission rates satisfy $\beta_m\in L^\infty_+(\R_+^2)$, $\beta_h\in L^\infty_+(\R_+^2)$;
    \item the rates $\gamma_h$ and $k_h$ belong to $L^\infty_+(\R_+^2)$;
    \item The initial condition is such that $S_{h,0}\in L^\infty_+(\R_+)$, $I_{h,0}\in L^\infty_+(\R_+^2)$, $R_{h,0}\in L^\infty_+(\R_+^2)$, $S_{m,0}\in L^\infty_+(\R_+)$, $I_{m,0}\in L^\infty(\R_+^2)$, with $\int_0^\infty S_{h,0}(a)\d a>0$ and $\int_0^\infty S_{m,0}(a)\d a>0$.
\end{enumerate}
\end{assumption}

To state our first main result, let us introduce the Banach space 
$$\X=L^1(\R_+)\times L^1(\R_+^2)\times L^1(\R_+^2)\times L^1(\R_+)\times L^1(\R_+^2),$$
endowed with the product norm. Its positive cone is defined by
$$\X_+=L^1_+(\R_+)\times L^1_+(\R_+^2)\times L^1_+(\R_+^2)\times L^1_+(\R_+)\times L^1_+(\R_+^2).$$
We also consider the space
$$\X_{\ep}=\{(S_h,I_h,R_h,S_m,I_m)\in \X: \|S_h\|_{L^1}+\|I_h\|_{L^1}+\|R_h\|_{L^1}\geq \ep\}.$$
The existence and uniqueness of bounded solutions of \eqref{Eq:Model} reads as follows.
\begin{theorem}[existence, uniqueness and boundedness]\label{Thm:existence}
Suppose that Assumption \ref{Assump} holds. Let
$$\overline{\ep}\in\left(0,\frac{\Lambda_h}{\|\mu_h\|_{L^\infty}+\|\nu_h\|_{L^\infty}+\|k_h\|_{L^\infty}+\|\gamma_h\|_{L^\infty}}\right).$$
Then Problem \eqref{Eq:Model} generates a globally defined strongly continuous semiflow $\{U(t)\}_{t\geq 0}$ on $\X_+\cap\X_{\overline{\ep}}$. For each $u_0=(S_{h,0}, I_{h,0}, R_{h,0}, S_{m,0}, I_{m,0})\in \X_+\cap \X_{\overline{\ep}}$, the total population of humans and mosquitoes at time $t$ satisfy the following inequalities:
\begin{equation}\label{Eq:Nh_gronw1}
  N_h(t)\leq N_h(0)e^{-\mu_0 t}+\frac{\Lambda_h}{\mu_0}\left(1-e^{-\mu_0 t}\right),
\end{equation}
\begin{equation}\label{Eq:Nh_gronw2}
N_h(t)\geq N_h(0)e^{-(\|\mu_h\|_{L^\infty}+\|\nu_h\|_{L^\infty})t}+\frac{\Lambda_h}{\|\mu_h\|_{L^\infty}+\|\nu_h\|_{L^\infty}}\left(1-e^{-(\|\mu_h\|_{L^\infty}+\|\nu_h\|_{L^\infty})t}\right),
\end{equation}
\begin{equation}\label{Eq:Nv_gronw1}
  N_m(t)\leq N_m(0)e^{-\mu_0 t}+\frac{\Lambda_m}{\mu_0}\left(1-e^{-\mu_0 t}\right),
\end{equation}
\begin{equation}\label{Eq:Nv_gronw2}
N_m(t)\geq N_m(0)e^{-(\|\mu_m\|_{L^\infty}+\|\nu_m\|_{L^\infty})t}+\frac{\Lambda_m}{\|\mu_h\|_{L^\infty}+\|\nu_m\|_{L^\infty}}\left(1-e^{-(\|\mu_m\|_{L^\infty}+\|\nu_m\|_{L^\infty})t}\right).
\end{equation}
Each population is bounded as follows:
\begin{equation}\label{Eq:bounded_human}
    \limsup_{t\to\infty}\|S_h(t,\cdot)\|_{L^1(\R_+)}\leq \frac{\Lambda_h}{\mu_0}, \quad \limsup_{t\to\infty}\|I_h(t,\cdot,\cdot)\|_{L^1(\R_+^2)}\leq \frac{\Lambda_h}{\mu_0}, \quad \limsup_{t\to\infty}\|R_h(t,\cdot,\cdot)\|_{L^1(\R_+^2)}\leq \frac{\Lambda_h}{\mu_0},
\end{equation}
\begin{equation}\label{Eq:bounded_mector}
\limsup_{t\to\infty}\|S_m(t,\cdot)\|_{L^1(\R_+)}\leq \frac{\Lambda_m}{\mu_0}, \quad \limsup_{t\to\infty}\|I_m(t,\cdot,\cdot)\|_{L^1(\R_+^2)}\leq \frac{\Lambda_m}{\mu_0}.
\end{equation}
Moreover, the following estimates hold:
\begin{equation}\label{Eq:bounded_Sh}
\liminf_{t\to \infty} S_h(t,a) \geq \Lambda_h \exp\left(-\int_0^a \left[\mu_h(s)+\frac{\Lambda_m\theta \|\beta_m\|_{L^\infty(\R_+^2)}(\|\mu_h\|_{L^\infty}+\|\nu_h\|_{L^\infty})}{\Lambda_h \mu_0}\right]\d s\right),
\end{equation}
\begin{equation}\label{Eq:bounded_Sv}
\liminf_{t\to \infty} S_m(t,a)\geq \Lambda_m\exp\left(-\int_0^a \left[\mu_m(s)+\theta\|\beta_h\|_{L^\infty}\right]\d s\right),
\end{equation}
for each $a\geq 0$. Finally, for each $t\geq 0$ we have
$$U(t)u_0=(S_h(t,\cdot), I_h(t,\cdot,\cdot), R_h(t,\cdot,\cdot), S_m(t,\cdot) , I_m(t,\cdot,\cdot)),$$
where components satisfy the following Volterra integral formulation:
\begin{equation*}
\begin{split}
S_h(t,a)= & \left\{
\begin{array}{lll}
   \int_0^a \left(\int_0^\infty k_h(s,\eta)R_h(t+s-a,s,\eta)\d \eta\right)e^{-\int_{s}^a \left(\mu_h(\xi)+\frac{\int_0^\infty \int_0^\infty \theta \beta_m(\zeta,\tau)I_m(t+s-\xi,\zeta,\tau)\d \zeta \d \tau}{N_h(t+s-\xi)}\right)\d \xi}\d s \\
   \quad + \Lambda_h e^{-\int_0^a \left(\mu_h(s)+\frac{\int_0^\infty \int_0^\infty \theta \beta_m(\xi,\tau)I_m(t+s-a,\xi,\tau)\d \xi \d \tau}{N_h(t+s-a)}\right)\d s}; \quad t>a,  \vspace{0.1cm} \\
 \int_0^t \left(\int_0^\infty k_h(a+s-t,\eta) R_h(s,a+s-t,\eta)\d \eta \right)e^{-\int_s^t \left(\mu_h(a+\xi-t)+\frac{\int_0^\infty \int_0^\infty 
\theta \beta_m(\zeta,\tau)I_m(\xi,\zeta,\tau)\d \zeta \d \tau}{N_h(\xi)}\right)\d \xi}\d s \\
    \quad + S_{h,0}(a-t) e^{-\int_{a-t}^a \left(\mu_h(s)+\frac{\int_0^\infty \int_0^\infty \theta \beta_m(\xi,\tau)I_m(t+s-a,\xi,\tau)\d \xi \d \tau}{N_h(t+s-a)}\right)\d s}; \quad a\geq t,
    \end{array}
    \right.\\
S_m(t,a)=&\left\{
\begin{array}{lll}
\displaystyle \Lambda_m e^{-\int_0^a \left(\mu_m(s)+\frac{\int_0^\infty \int_0^\infty \theta \beta_h(\xi,\tau)I_h(t+s-a,\xi,\tau)\d \xi \d \tau}{N_h(t+s-a)}\right)\d s}, \quad t>a,  \vspace{0.1cm} \\
\displaystyle S_{m,0}(a-t)e^{-\int_{a-t}^a \left(\mu_m(s)+\frac{\int_0^\infty \int_0^\infty \theta \beta_h(\xi,\tau)I_h(t+s-a,\xi,\tau)\d \xi \d \tau}{N_h(t+s-a)}\right)\d s}; \quad a\geq t,
\end{array}
\right.   
\end{split}
\end{equation*}

\begin{equation*}
\begin{split}
I_h(t,a,\tau)= & \left\{
\begin{array}{lll}
  \left(\frac{S_h(t-\tau,a-\tau)}{N_h(t-\tau)}e^{-\int_0^\tau (\mu_h(s+a-\tau)+\nu_h(s+a-\tau,s)+\gamma_h(s+a-\tau,s))\d s}\right) \\
  \qquad \times \displaystyle\left(\int_0^\infty\int_0^\infty \theta \beta_m(s,\xi)I_m(t-\tau,s,\xi)\d s\d \xi\right); \quad t>\tau; \quad a\geq \tau,  \vspace{0.1cm} \\
  I_{h,0}(a-t,\tau-t)e^{-\int_{\tau-t}^\tau (\mu_h(s+a-\tau)+\nu_h(s+a-\tau,s)+\gamma_h(s+a-\tau,s))\d s}; \quad a\geq \tau\geq t,
\end{array}
\right.\\
I_m(t,a,\tau)=& \left\{
\begin{array}{lll}
 \left(\frac{S_m(t-\tau,a-\tau)}{N_h(t-\tau)}\int_0^\infty \int_0^\infty \theta \beta_h(s,\xi)I_h(t-\tau,s,\xi)\d s~\d \xi\right)\\
 \qquad \times \displaystyle \left(e^{-\int_0^\tau (\mu_m(s+a-\tau)+\nu_m(s+a-\tau,s))\d s}\right); \ t>\tau; \quad a\geq \tau,  \vspace{0.1cm} \\
 I_{m,0}(a-t,\tau-t)e^{-\int_{\tau-t}^\tau (\mu_m(s+a-\tau)+\nu_m(s+a-\tau,s))\d s}; \ a\geq \tau\geq t,
\end{array}
\right. \\
R_h(t,a,\eta)=& \left\{
\begin{array}{lll}
 \left(\int_0^\infty \gamma_h(a-\eta,\tau)I_h(t-\eta,a-\eta,\tau)\d \tau\right)e^{-\int_0^\eta (\mu_h(s+a-\eta)+k_h(s+a-\eta,s))\d s}; \ t>\eta; \ a\geq \eta,  \vspace{0.1cm} \\
 \displaystyle R_{h,0}(a-t,\eta-t)e^{-\int_{\eta-t}^\eta (\mu_h(s+a-\eta)+k_h(s+a-\eta,s))\d s}; \quad a\geq \eta\geq t.
\end{array}
\right.
\end{split}
\end{equation*}
\end{theorem}

\subsection{The disease invasion process}
We see that there always exists a disease-free equilibrium denoted by
$$E_0=\left(S^0_h, I^0_h=0, R^0_h=0, S^0_m, I^0_m=0\right)$$
where
$$S^0_h(a)=\Lambda_h e^{-\int_0^a \mu_h(s)\d s}, \qquad S^0_m(a)=\Lambda_m e^{-\int_0^a \mu_m(s)\d s} \qquad \forall a\in\R_+.$$ Note that $$\pi_h(a)= e^{-\int_0^a \mu_h(s)\d s}, \qquad \pi_m(a)= e^{-\int_0^a \mu_m(s)\d s}$$ are survival rates, from birth until age $a$, for humans and mosquitoes respectively, in absence of disease. Therefore, $S^0_h(a)$ and $S^0_m(a)$ are average number of humans and mosquitoes aged $a$ in a disease-free environment. 

The number of new infections in humans that one human causes through his/her infectious period is given by $\mathcal{R}_0^2$, where $\mathcal{R}_0$ is the basic reproduction number characterized as the spectral radius  based on the next generation operator approach \cite{Diekmann90,Inaba2012}. However, we set the following technical assumption:
\begin{assumption}\label{Assum:beta}
we suppose that $\beta_h$ and $\beta_m$ are not identically zero on the set $\{(s+\tau,\tau): (s,\tau)\in \R_+^2\}$.
\end{assumption}
Note that biologically the age is always larger than the time since infection. Consequently, the latter assumption only implies that infected humans and mosquitoes will be infectious at some point of the infection.

We find that, under Assumptions \ref{Assump} and \ref{Assum:beta}, the basic reproduction number is given by
\begin{flalign}
\mathcal{R}_0^2= & \underbrace{\frac{\Lambda_m \int_0^\infty \pi_m(s)\d s}{\Lambda_h \int_0^\infty \pi_h(s)\d s}}_{\text{mosquito/human ratio}} \times \underbrace{\theta^2  \int_0^\infty \int_0^\infty \mathcal{K}_{m\to h}(\xi,\tau)   \d \xi~\d \tau}_{\text{transmission rate from mosquito to human}} \nonumber   \times \underbrace{ \int_0^\infty \int_0^\infty  \mathcal{K}_{h\to m}(\xi,\tau) \d \xi~\d \tau}_{\text{per bite transmission rate from human to mosquito}} \label{Eq:R0}
\end{flalign}
where
\begin{align*}
    \mathcal{K}_{h\to m}(\xi,\tau)=& \beta_h(\tau+\xi,\tau) e^{-\int_0^\tau (\nu_h(\sigma+\xi,\sigma)+\gamma_h(\sigma+\xi,\sigma))\d \sigma} \frac{\pi_h(\xi+\tau)}{\int_0^\infty \pi_h(s)\d s},\\
    \mathcal{K}_{m \to h}(\xi,\tau)=& \beta_m(\tau+\xi,\tau)  e^{-\int_0^\tau \nu_m(\sigma+\xi,\sigma)\d \sigma} \frac{\pi_m(\xi+\tau)}{\int_0^\infty \pi_m(s)\d s}.
\end{align*}

Here, $\mathcal{K}_{h\to m}(\xi,\tau)$ describes the infectiousness of a human with age $\xi+ \tau$ and time since infection $\tau$. It quantifies the proportion of the bites by susceptible mosquitoes on infectious humans, with age $\xi+ \tau$ and time since infection $\tau$, that infect mosquitoes. More precisely, it is given by the product between $\beta_h(\tau+\xi,\tau)$ the disease transmission probability, and
$$e^{-\int_0^\tau (\nu_h(\sigma+\xi,\sigma)+\gamma_h(\sigma+\xi,\sigma))\d \sigma} \frac{\pi_h(\xi+\tau)}{\int_0^\infty \pi_h(s)\d s}$$
that is the survival rate of an infected human with age $\xi+ \tau$ and time since infection $\tau$. Similarly, $\mathcal{K}_{m\to h}(\xi,\tau)$ is the infectiousness of a mosquito with age $\xi+ \tau$ and time since infection $\tau$, {\it i.e.,} the proportion of the bites by infectious mosquitoes, with age $\xi+ \tau$ and time since infection $\tau$, that infect susceptible humans. Once multiplying $\mathcal{K}_{m\to h}(\xi,\tau)$ by mosquitoes biting rate and integrating over all chronological and infection ages $\xi$ and $\tau$ it gives the vectorial capacity (or the the ability of the vector to transmit the disease).  

From above notations, we now state our disease invasion threshold criterion as follows
\begin{theorem}[stability of $E_0$]\label{Thm:stab}
Under Assumptions \ref{Assump} and \ref{Assum:beta} we have:
\begin{enumerate}
    \item if $\mathcal{R}_0>1$, then $E_0$ is unstable;
    \item if $\mathcal{R}_0<1$, then $E_0$ is locally asymptotically stable.
\end{enumerate}
\end{theorem}

\subsection{Existence and bifurcation of an endemic equilibrium}
Let $$E^*=\left(S^*_h(a),I^*_h(a,\tau), R^*_h(a,\eta), S^*_m(a), I^*_m(a,\tau)\right)>0$$ be an endemic equilibrium of \eqref{Eq:Model}. Setting $$N^*_h= \int_0^\infty S^*_h(a)\d a+\int_0^\infty \int_0^\infty I^*_h(a,\tau)\d a~\d \tau+\int_0^\infty \int_0^\infty R^*_h(a,\eta)\d a~\d \eta$$ and 
$$s^*_h=\frac{S^*_h}{N^*_h}, \qquad i^*_h=\frac{I^*_h}{N^*_h}, \qquad r^*_h=\frac{R^*_h}{N^*_h}$$ we find that a necessary and sufficient condition for the existence of $E^*$ is given by 
\begin{equation*}
\begin{array}{ll}
1=& \mathcal{R}_0^2\left(1+\int_0^\infty \int_0^a \left(\int_0^s \nu_h(s,\tau) i_h^*(s,\tau)\d \tau\right) \exp\left(-\int_s^a \mu_h(\xi)\d \xi\right)\d s~\d a\right)^2\\
&\times\left(\frac{\int_0^\infty \int_\tau^\infty \beta_m(a,\tau)\exp\left(-\int_0^a \mu_m(s)\d s-\int_0^\tau \nu_m(s+a-\tau,s)\d s+(\tau-a)\int_0^\infty \int_\xi^\infty \theta \beta_h(s,\xi)i^*_h(s,\xi)\d s~\d \xi\right)\d a~\d \tau}{\int_0^\infty \int_\tau^\infty \beta_m(a,\tau)\exp\left(-\int_0^a \mu_m(s)\d s-\int_0^\tau \nu_m(s+a-\tau,s)\d s\right)}\right)  \\
&-\left(\int_0^\infty \int_\tau^\infty \theta \beta_m(a,\tau)e^{-\int_0^a \mu_m(s)\d s-\int_0^\tau \nu_m(s+a-\tau,s)\d s+(\tau-a)\int_0^\infty \int_\xi^\infty \theta \beta_h(s,\xi)i^*_h(s,\xi)\d s~\d \xi}\d a~\d \tau\right) \\
&\times \Lambda_m\left(\frac{1+\int_0^\infty \int_0^a \left(\int_0^s \nu_h(s,\tau) i_h^*(s,\tau)\d \tau\right) \exp\left(-\int_s^a \mu_h(\xi)\d \xi\right)\d s~\d a}{\Lambda_h \int_0^\infty \exp\left(-\int_0^a \mu_h(s)\d s\right)\d a}\right) \\
&\times \int_0^\infty \int_{\tau}^\infty \theta \beta_h(a,\tau)\exp\left(-\int_0^\tau (\mu_h(s+a-\tau)+\nu_h(s+a-\tau,s)+\gamma_h(s+a-\tau,s))\d s\right)\\
&\times \left[\int_0^{a-\tau} i_h^*(a-\tau,s)\d s+\int_0^{a-\tau}\left(\int_0^s \nu_h(s,\tau)i^*_h(s,\tau)\d \tau\right)e^{-\int_s^{a-\tau} \mu_h(\xi) \d \xi}\d s\right. \\
&\left.+\int_0^{a-\tau} \left(\int_0^\infty \gamma_h(a-\tau-\eta,\tau)i^*_h(a-\tau-\eta,\tau)\d \tau\right)e^{-\int_0^\eta (\mu_h(s+a-\tau-\eta)+k_h(s+a-\tau-\eta,s))\d s} \right]\d a~\d \tau .
\end{array}
\end{equation*}
The above condition has only $i_h^*$ as unknown, but is though difficult to exploit, particularly because $i_h^*$ depends both on $a$ and $\tau$. To go further, we additionally assume that: 
\begin{assumption}\label{Assump-relax-human-age}
$\mu_h, \nu_h, \gamma_h, k_h$ and $\beta_h$ do not depend on humans age, {\it i.e.,} $\mu_h(a) \equiv \mu_h$, $\nu_h(a,\tau) \equiv \nu_h(\tau)$, $\gamma_h(a,\tau) \equiv \gamma_h(\tau)$, $k_h(a,\eta)=k_h(\eta)$ and $\beta_h(a,\tau) \equiv \beta_h(\tau)$.
\end{assumption}

First, let us note that under Assumption \ref{Assump-relax-human-age}, the previous necessary and sufficient condition for the existence of $E^*$ simplified. Indeed, we find that an endemic equilibrium of \eqref{Eq:Model} exists if and only if there exists $K>0$ such that
$$f(\mathcal{R}_0^2,K)=1$$
with $f$ defined by 
\begin{flalign}
&f(R,K)=R\left(1+\frac{K\int_0^\infty \nu_h(\tau)\pi_h(\tau) e^{-\int_0^\tau (\nu_h(s)+\gamma_h(s))\d s}\d \tau}{\mu_h}\right)\label{Eq:function_f} \\
&\left(\int_0^\infty \int_\tau^\infty  \beta_m(a,\tau)\pi_m(a) e^{-\int_0^\tau \nu_m(s+a-\tau,s)\d s} e^{-(a-\tau)K\int_0^\infty \theta \beta_h(s)\pi_h(s) e^{-\int_0^s (\nu_h(\xi)+\gamma_h(\xi))\d \xi}\d s}\d a~\d \tau\right) \nonumber \\
&\times \left[\frac{1-K\int_0^\infty \pi_h(\tau) e^{-\int_0^\tau (\nu_h(s)+\gamma_h(s))\d s}\left(1+\gamma_h(\tau)\int_0^\infty \pi_h(\eta) e^{-\int_0^\eta k_h(s)\d s}\d \eta\right)\d \tau}{\int_0^\infty \int_\tau^\infty \beta_m(a,\tau)\pi_m(a) e^{-\int_0^\tau \nu_m(s+a-\tau,s)\d s} \d a~\d \tau}\right]. \nonumber
\end{flalign}

Next, let us also introduce the following bifurcation constant
\begin{flalign}
C_{\bif}=&\frac{\int_0^\infty \beta_h(\tau) \pi_h(\tau) e^{-\int_0^\tau(\nu_h(s)+\gamma_h(s))\d s}\d \tau}{\int_0^\infty \int_\tau^\infty \beta_m(a,\tau) \pi_m(a) e^{-\int_0^\tau \nu_m(s+a-\tau,s)\d s}\d a~\d \tau} \int_0^\infty \int_\tau^\infty \theta \beta_m(a,\tau) (\tau-a) \pi_m(a) e^{-\int_0^\tau \nu_m(s+a-\tau,s)\d s}\d a~\d \tau \nonumber \\
&-\int_0^\infty \gamma_h(\tau)e^{-\int_0^\tau(\mu_h+k_h(s))\d s}\d \tau +\int_0^\infty e^{-\int_0^\tau \pi_h(\tau) (\nu_h(s)+\gamma_h(s))\d s}\left(\frac{\nu_h(\tau)}{\mu_h}-1\right)\d \tau. \label{Eq:constant_bif}
\end{flalign}

We then have the following existence and bifurcation result of endemic equilibrium
\begin{theorem}[Bifurcations]\label{Thm:Bif}
Let Assumptions \ref{Assump} and \ref{Assump-relax-human-age} be satisfied. It comes:
\begin{enumerate}
    \item if $C_\bif>0$, then there is a backward bifurcation at $\mathcal{R}_0=1$, \textit{i.e.} for $\mathcal{R}_0<1$ close enough to $1$, there exists two endemic equilibria;
    \item if $C_\bif<0$, then there is a forward bifurcation at $\mathcal{R}_0=1$, \textit{i.e.} for $\mathcal{R}_0>1$ close enough to $1$, there exists a unique endemic equilibrium and for $\mathcal{R}_0<1$ close enough to $1$, there is no endemic equilibrium;
     \item if $\mathcal{R}_0>1$ then there exists at least one endemic equilibrium of \eqref{Eq:Model}. Moreover, if $C_\bif>0$, then whenever $\mathcal{R}_0=1$, there is also an endemic equilibrium; 
    \item suppose that the following condition holds:
    \begin{equation}\label{Eq:particular_case}
    \int_0^\infty \pi_h(\tau) \exp\left(-\int_0^\tau (\nu_h(s)+\gamma_h(s))\d s\right)\left(\frac{\nu_h(\tau)}{\mu_h}-1\right)\d \tau\leq 0
    \end{equation}
    then there exists an endemic equilibrium if and only if $\mathcal{R}_0>1$, and in that case the equilibrium is unique.
\end{enumerate}
\end{theorem}

Note that when the condition \eqref{Eq:particular_case} is satisfied, we systematically have $C_\bif<0$. For example, the condition \eqref{Eq:particular_case} holds if the humans death rate induced by the infection is small enough, {\it i.e.} $\sup_{\tau} \nu_h(\tau) \leq \mu_h$.

Details on the proof of our main results are given after some numerical simulations to qualitatively illustrate such results.

\section{Numerical simulations} \label{Sec:simus}

\begin{figure}[!h]
\centering
\begin{tabular}{cc}
\includegraphics[width=.45\linewidth]{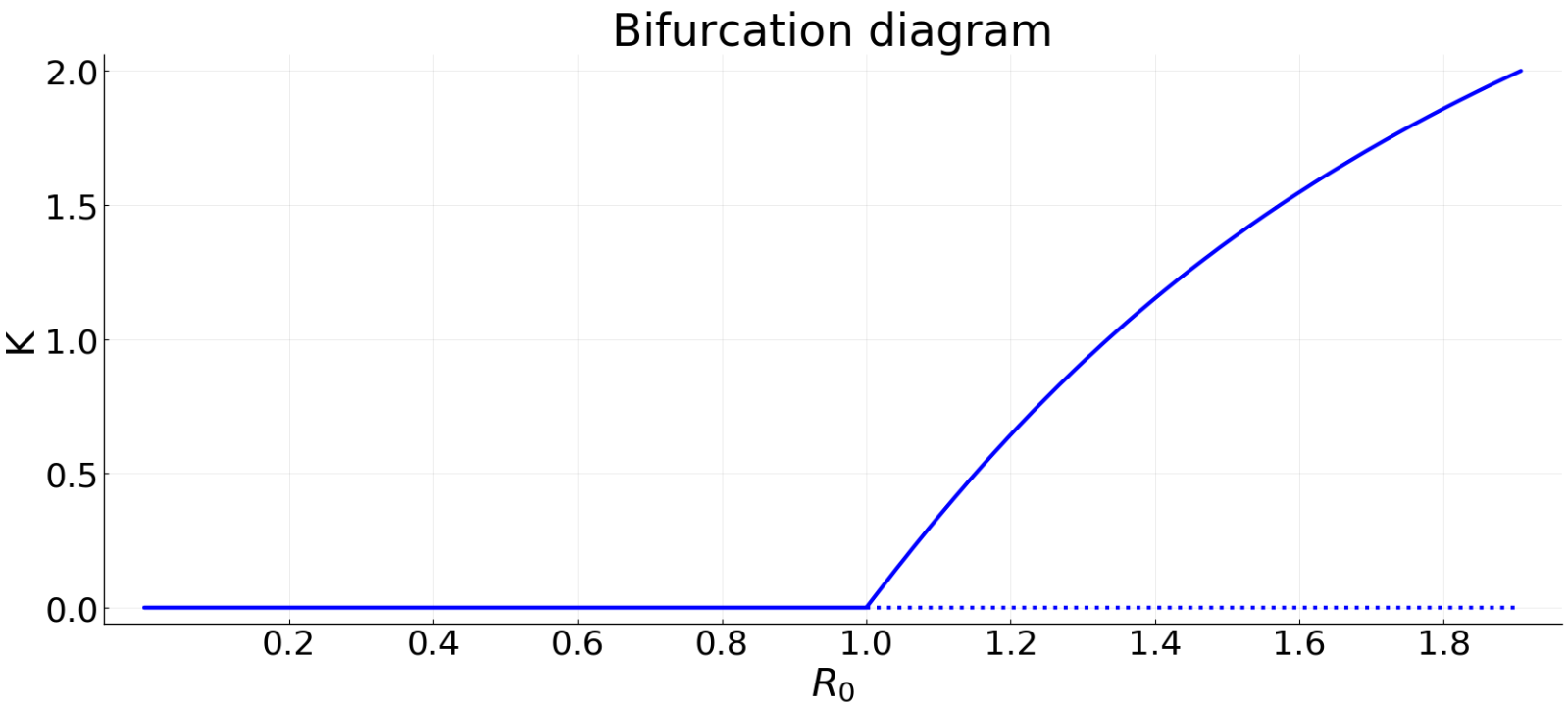} & \includegraphics[width=.45\linewidth]{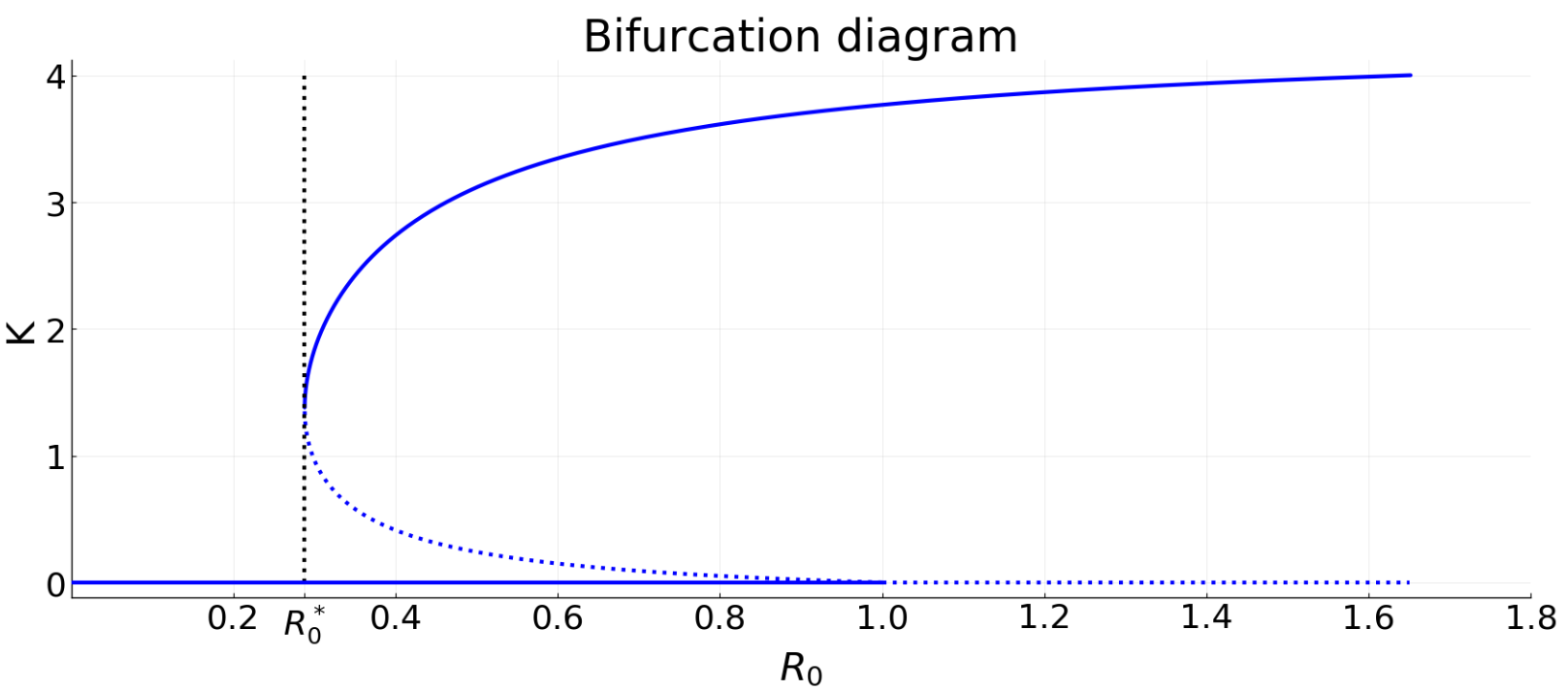} 
\end{tabular}
\caption{Forward bifurcation when $\mu_h=0.022$ (left) and backward bifurcation when $\mu_h=0.002$ (right).} \label{Fig:bif}
\end{figure}

In this section we show some numerical simulations, by using finite volume numerical schemes (implemented with the Julia Programming Language) to illustrate bifurcation results of the endemic equilibrium of Model \eqref{Eq:Model} under Assumptions \ref{Assump}-\ref{Assump-relax-human-age}. 
We randomly set the parameters, with the only purpose to illustrate the bifurcations results. Note that the parametrisation of the model, with existing data on the rates will be addressed in a further work. First we fix the following parameters:
$$\Lambda_h=8.4\times 10^5, \quad \mu_m(a)=20, \quad \nu_h(a)=0.1, \quad \nu_m(a,\tau)=25, \quad \theta=3.65 \times 10^4,$$
for each $(a,\tau)\in \R_+^2$,
$$\gamma_h(\tau)=
\begin{cases}
0 & \text{if } \tau\in[0,0.1],
\\ 50 &\text{else }
\end{cases}, \qquad k_h(\eta)=
\begin{cases}
0 & \text{if } \eta \in[0,0.1],
\\ 40 &\text{else }
\end{cases} \qquad \beta_h(\tau)= 	\frac{0.1}{\sqrt{2\pi}}e^{-\frac{1}{2}\left(\frac{\tau-0.3}{0.1}\right)^2} 
$$
and
$$\beta_m(a,\tau)=
\begin{cases}
0 & \text{if } a\leq \tau,
\\ \frac{0.05}{\sqrt{2\pi}}e^{-\frac{1}{2}\left(\frac{\tau-0.2}{0.2}\right)^2} e^{-(a-\tau)} & \text{ else}.
\end{cases}
$$
Now, we may observe that the bifurcation constant $C_{\bif}$ defined by \eqref{Eq:constant_bif}, depends on $\mu_h$ but not on $\Lambda_m$. We thus consider two cases:
\begin{enumerate}
    \item $\mu_h=0.022$ which leads to $C_{\bif}\approx -1.34<0$ and a forward bifurcation occurs at $\mathcal{R}_0=1$ according to Theorem \ref{Thm:Bif} (2);
    \item $\mu_h=0.002$ which leads to $C_{\bif}\approx 4>0$ and a backward bifurcation occurs at $\mathcal{R}_0=1$ according to Theorem \ref{Thm:Bif} (1).
\end{enumerate}
Moreover, the threshold $R_0:= \mathcal{R}_0^2$ can be written as a linear function of $\Lambda_m$, which will be considered as bifurcation parameter. Then, injecting it into the necessary and sufficient condition $f(R_0,K)=1$ (ensuring the existence of an endemic equilibrium) allows us to draw the bifurcation figures (see Figure \ref{Fig:bif}). We remind here that $K>0$ necessary implies, by definition, that $I^*_m\not\equiv 0$.

\subsection{Forward bifurcation}

\begin{figure}[!h]
\centering
\begin{tabular}{cc}
\includegraphics[width=.45\linewidth]{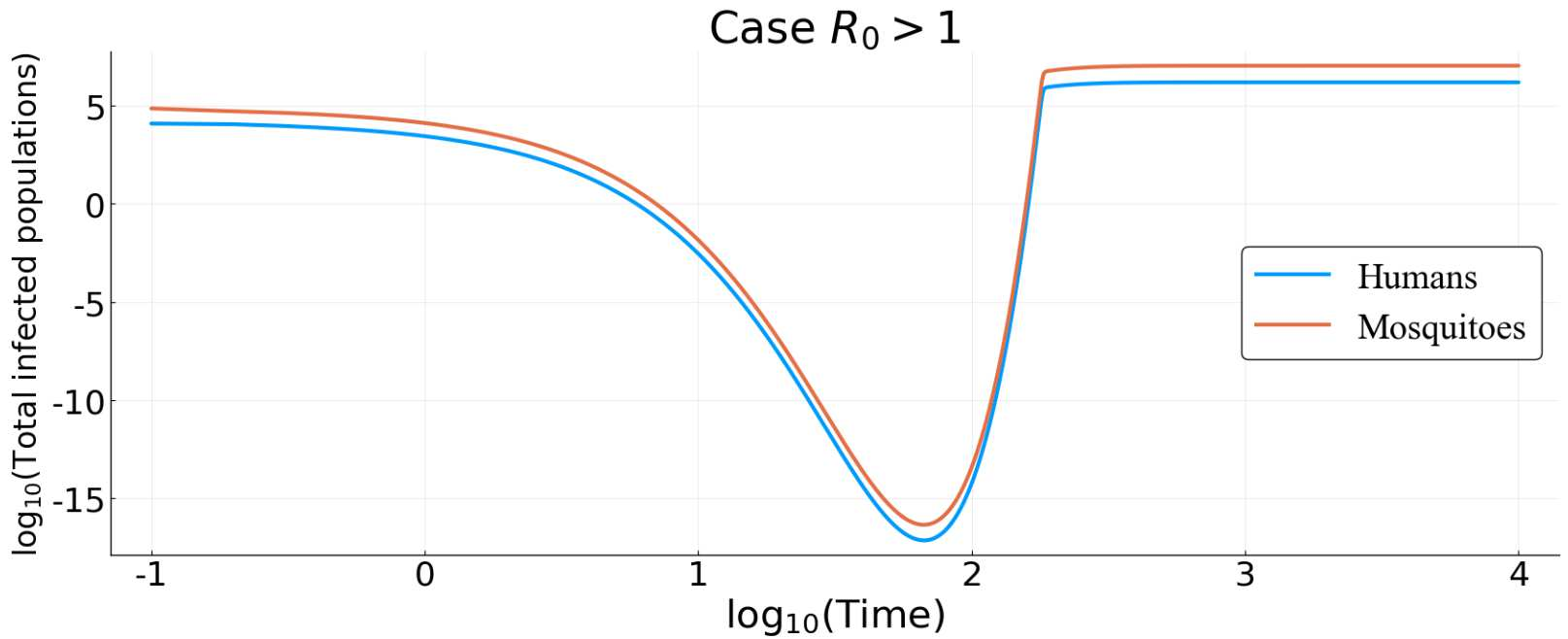} & \includegraphics[width=.45\linewidth]{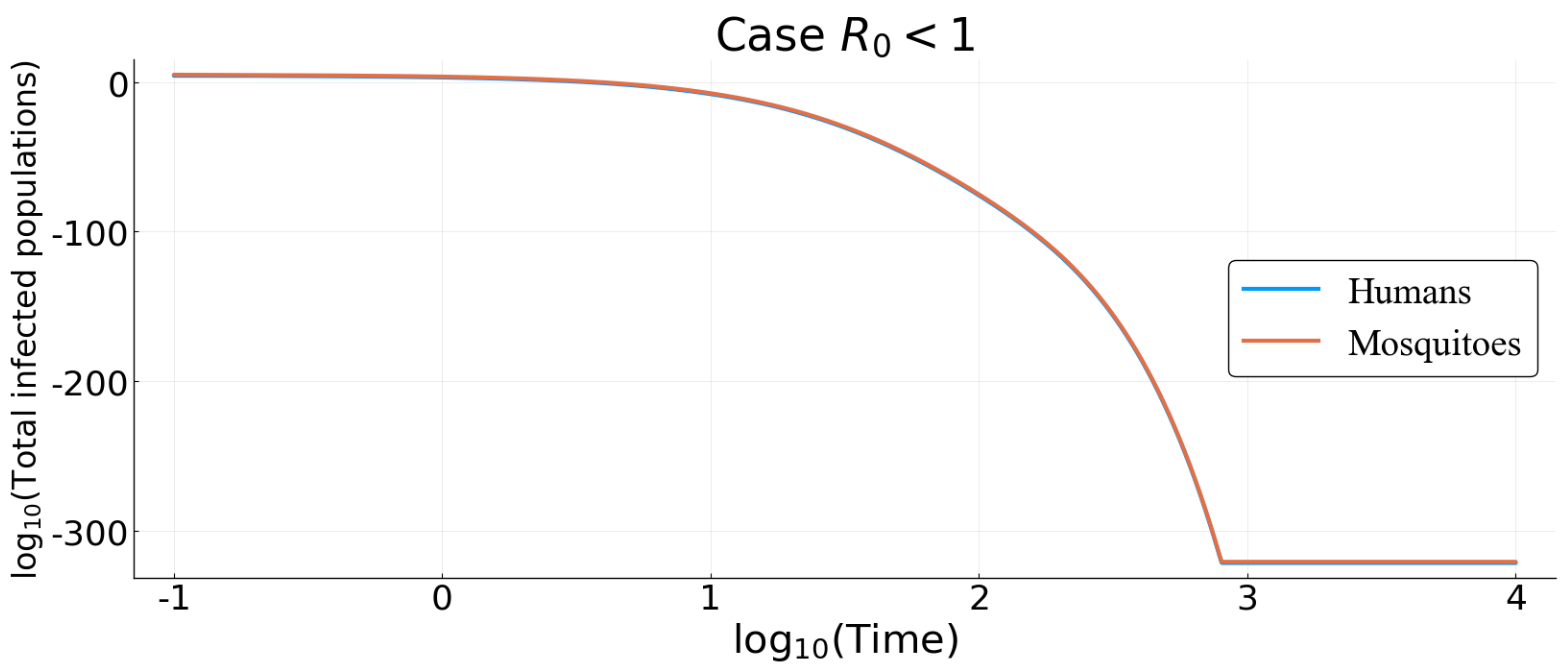} 
\end{tabular}
\caption{Case of forward bifurcation: either convergence to the endemic equilibrium when $R_0>1$ (left) or to the disease-free equilibrium when $R_0<1$ (right).} \label{Fig:convforw}
\end{figure}

Let us fix $\mu_h=0.022$. As mentioned above, a forward bifurcation occurs at $R_0=1$, which means that whenever $R_0<1$, then the disease-free equilibrium $E_0$ is locally asymptotically stable by Theorem \ref{Thm:stab} and no endemic equilibrium exists. Asymptotically, the disease will go extinct (see Figure \ref{Fig:convforw} left, where $\Lambda_m=7\times 10^6$ and $R_0\approx 1.16$). However, if $R_0>1$, then $E_0$ is unstable by Theorem \ref{Thm:stab} and an endemic equilibrium exists by Theorem \ref{Thm:Bif}. The disease will asymptotically persist and under Assumption \ref{Assump-relax-human-age} the solution of \eqref{Eq:Model} will converge to the endemic equilibrium $E^*$ (see Figure \ref{Fig:convforw} right, where $\Lambda_m=5\times 10^6$ and $R_0\approx 0.83$). Note that in Figure \ref{Fig:convforw} (as well as Figure \ref{Fig:convback}), both axes are in $\log_{10}$ scale with time in $x$-axis and the total population of infected in $y$-axis for both humans and mosquitoes.

\subsection{Backward bifurcation}

\begin{figure}[!h]
\centering
\begin{tabular}{cc}
\includegraphics[width=.45\linewidth]{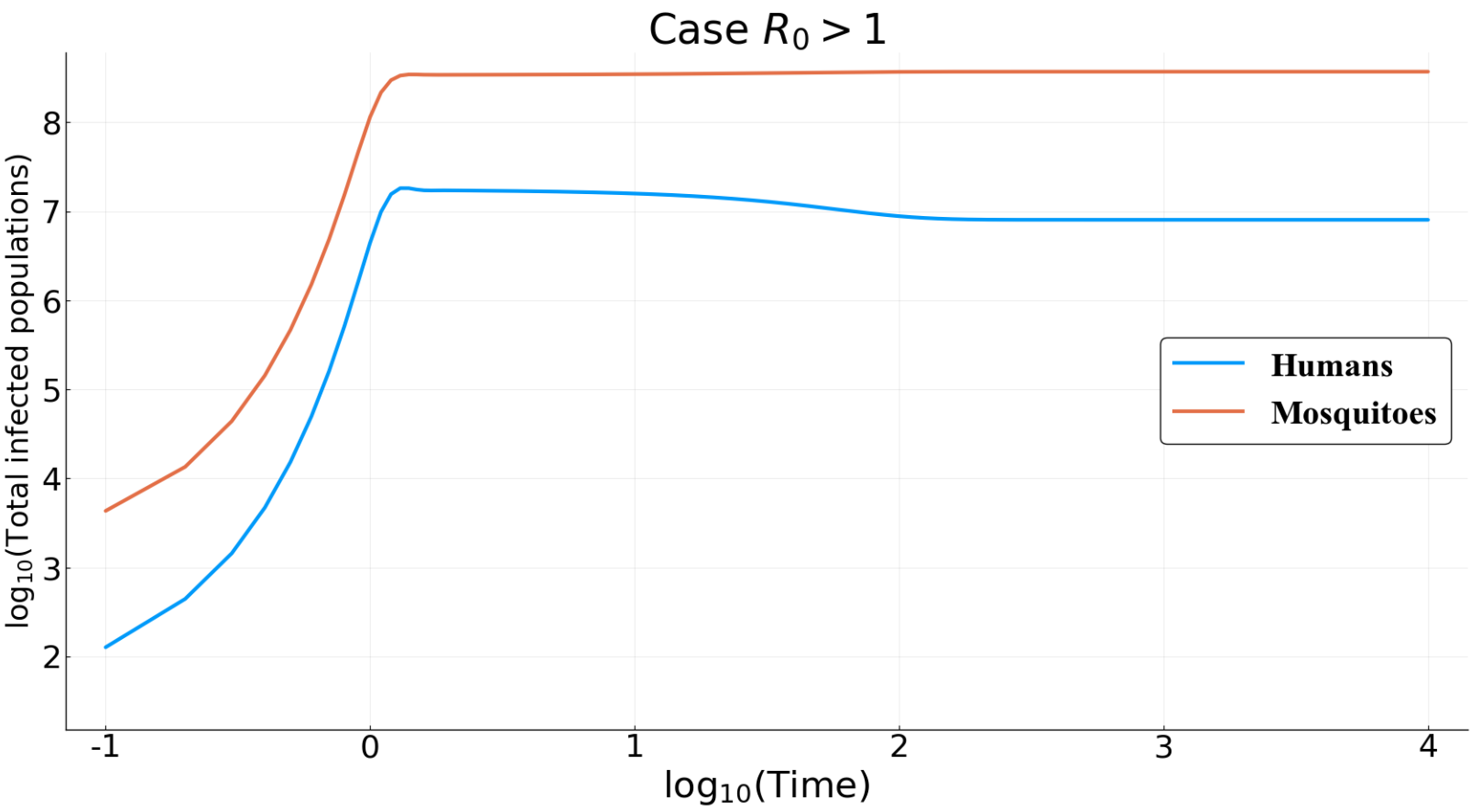} & \includegraphics[width=.45\linewidth]{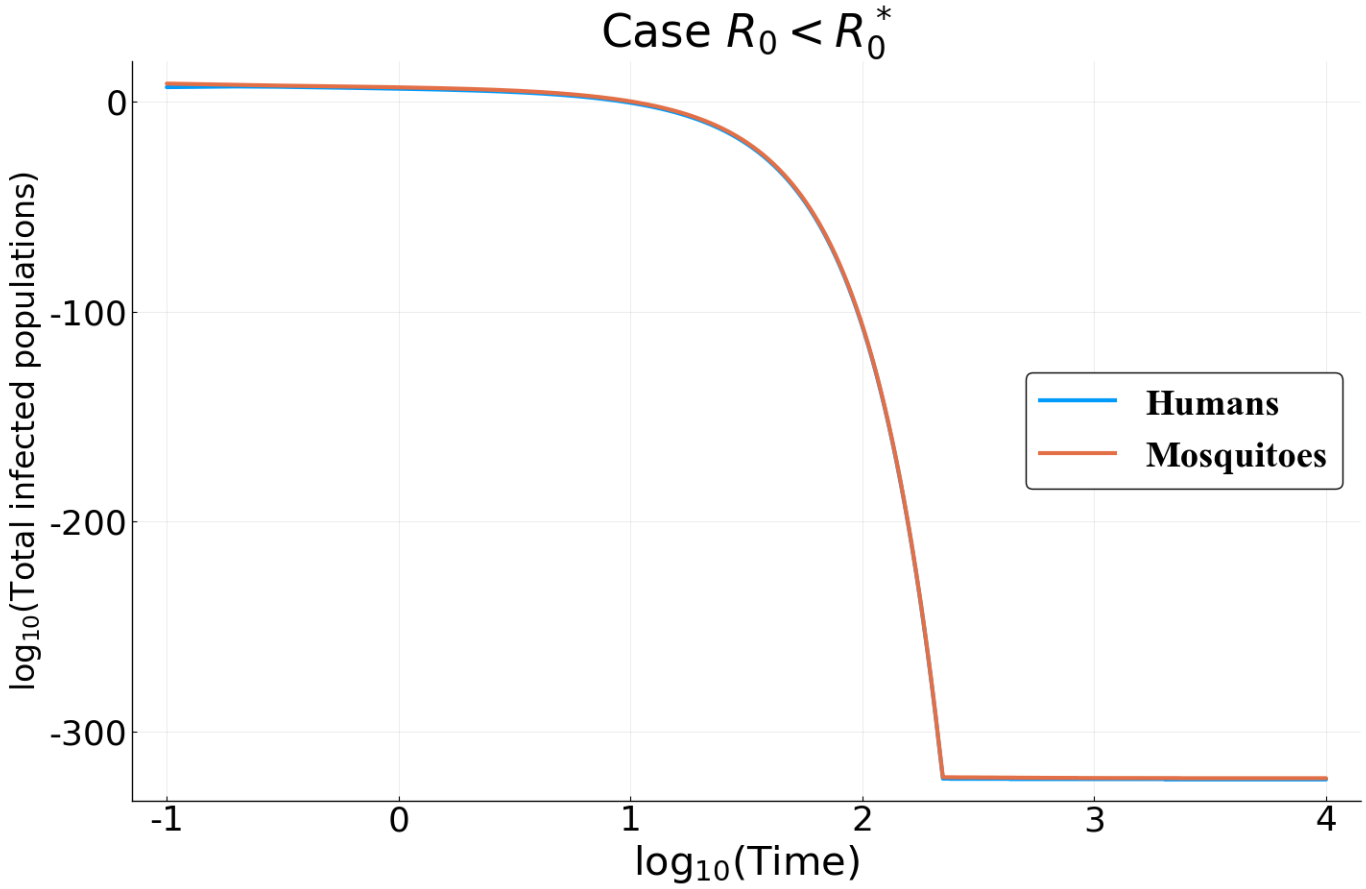}  \\
\includegraphics[width=.45\linewidth]{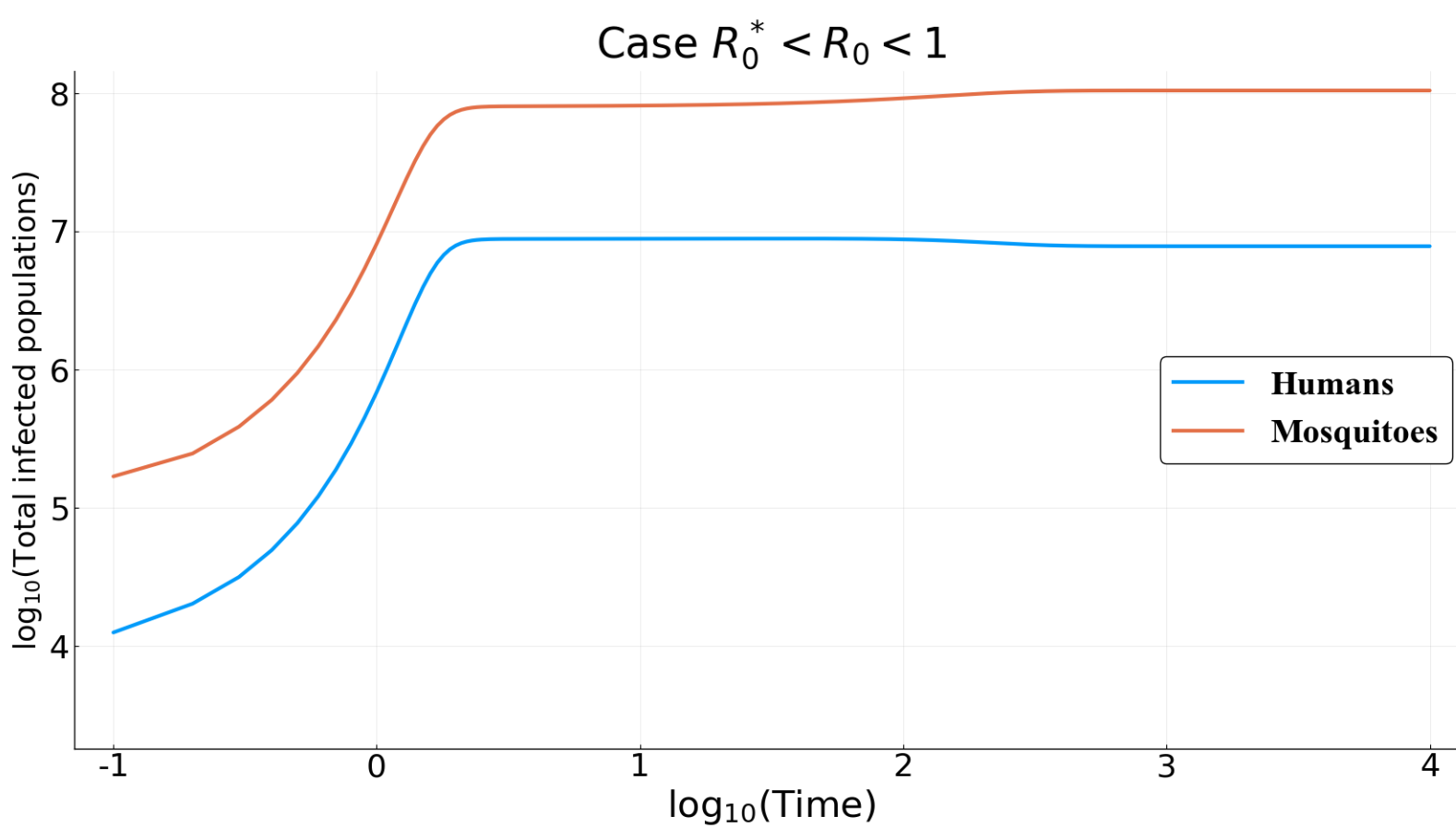} & \includegraphics[width=.45\linewidth]{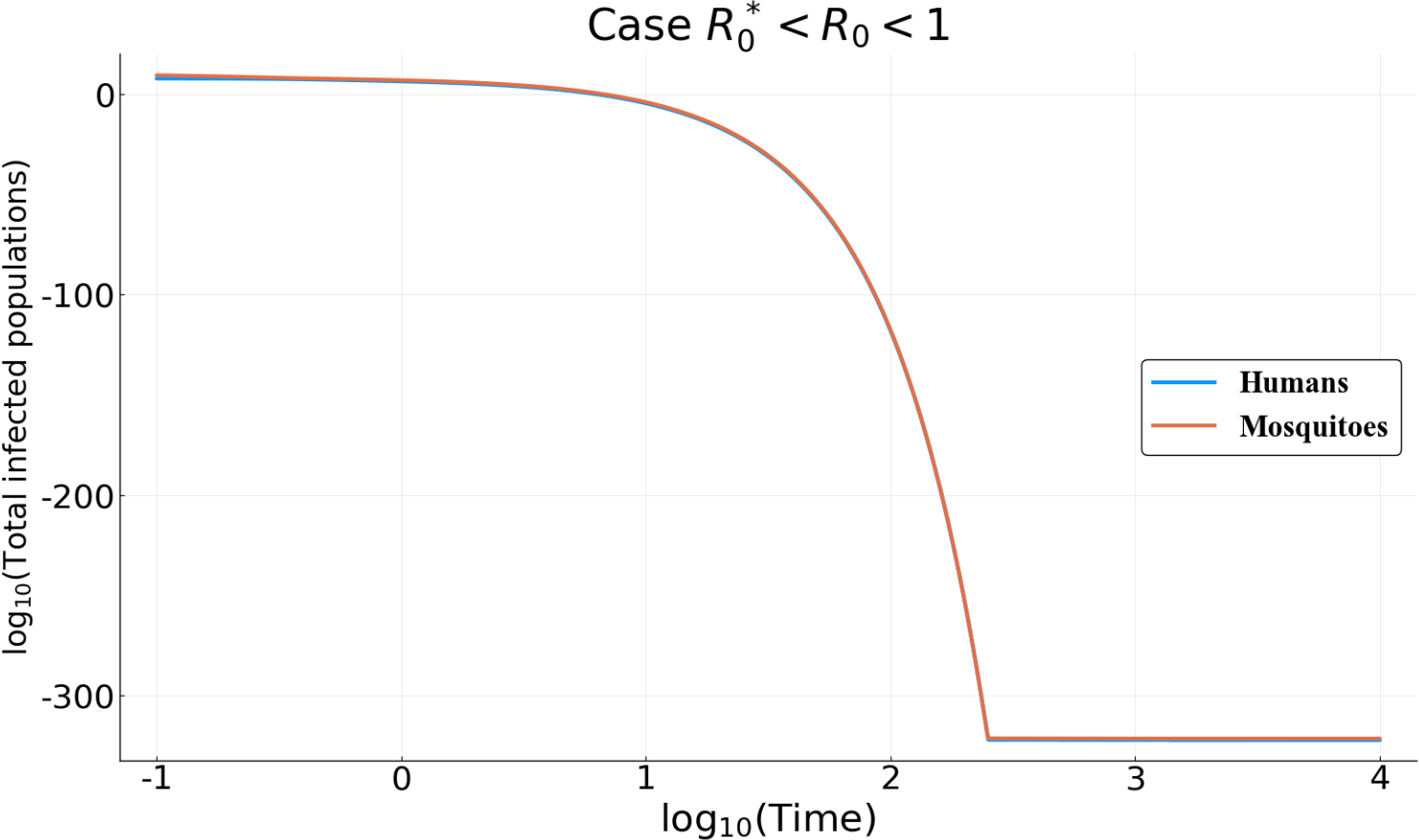} 
\end{tabular}
\caption{Case of backward bifurcation: either convergence to the endemic equilibrium when $R_0>1$ (top left) or to the disease-free equilibrium when $R_0<R_0^*<1$ (top right). In the case $R_0^*<R_0<1$ and considering two different initial conditions: either convergence to an endemic equilibrium (bottom left) or to the disease-free equilibrium (bottom right).} \label{Fig:convback}
\end{figure}

Here we fix $\mu_h=0.002$, where a backward bifurcation occurs at $R_0=1$. It implies that whenever $R_0>1$, the disease-free equilibrium is unstable and there exists a unique endemic equilibrium. Asymptotically, the solutions converge to this endemic equilibrium (see Figure \ref{Fig:convback} top left, where $\Lambda_m=7.4 \times 10^7$ and $R_0\approx 1.12$). However, contrary to the forward bifurcation, there exists a threshold $R_0^*<1$ (numerically $R_0^*\approx 0.29$) such that: if $R_0<R_0^*$, then there is no endemic equilibrium and under Assumption \ref{Assump-relax-human-age}, the solutions of \eqref{Eq:Model} converge to the disease-free equilibrium $E_0$ (see Figure \ref{Fig:convback} top right, where $\Lambda_m=10^7$ and $R_0\approx 0.15$). 

When $R_0^*<R_0<1$, then there exist at least two endemic equilibria with the corresponding two values of $K$. We set $\Lambda_m=2.5\times 10^7$, so that $R_0\approx 0.38\in[R_0^*,1]$. Considering two sets of initial conditions for $(I_{h,0}, I_{m,0})$ and under Assumption \ref{Assump-relax-human-age}, the solutions of \eqref{Eq:Model} will either converge to one of the endemic equilibria (see Figure \ref{Fig:convback} bottom left), or to the disease-free equilibrium $E_0$ (see Figure \ref{Fig:convback} bottom right).

\section{Proof of Theorem \ref{Thm:existence}: existence and uniqueness of bounded solutions}\label{Sec:existence}
In this section, under Assumption \ref{Assump}, we handle the well-posedness of \eqref{Eq:Model} with an integrated semigroups approach.
\subsection{Integrated semigroup formulation}\label{Sec:Sg_formulation}

We introduce the space
$$X_1=\R\times L^1(\R_+),$$
and the linear operators $A_{h,1}:D(A_{h,1})\subset X_1\to X_1$ and $A_{m,1}:D(A_{m,1})\subset X_1\to X_1$ be defined by
$$A_{h,1}\begin{pmatrix}
0_{\R} \\
\phi
\end{pmatrix}=\begin{pmatrix}
-\phi(0)\\
-\phi'-\mu_h \phi\end{pmatrix}, \qquad A_{m,1}\begin{pmatrix}
0_{\R} \\
\phi
\end{pmatrix}=\begin{pmatrix}
-\phi(0)\\
-\phi'-\mu_m \phi\end{pmatrix},
$$
where the domain $D(A_{h,1})$, $D(A_{m,1})$ of operators $A_{h,1}$, $A_{m,1}$ are $D(A_{h,1})=D(A_{m,1})=\{0_{\R}\}\times W^{1,1}(\R_+)$. 

By Assumption \ref{Assump} on mortality rates, we find that if $\lambda\in \C$ is such that $\Re(\lambda)\geq -\mu_0$, then $\lambda\in \rho(A_{h,1})\cap \rho(A_{m,1})$ (with $ \rho(A)$ the resolvent of any operator $A$), and we have the following explicit formula for the resolvent of $A_{k,1}$ (with $k\in\{h,m\})$:
\begin{equation}\label{Eq:Resolv1}
    \left(\lambda I-A_{k,1}\right)^{-1}\begin{pmatrix}
    c \\
    \psi
    \end{pmatrix}
    =\begin{pmatrix}
    0 \\
    \phi
    \end{pmatrix}\Longleftrightarrow \phi(a)=c e^{-\int_0^a (\mu_k(s)+\lambda)\d s}+\int_0^a \psi(s)e^{-\int_s^a (\mu_k(\xi)+\lambda)\d \xi}\d s,
\end{equation}
for $(c,\psi)^T\in X$. Now, we introduce the space
$$X_2=L^1(\R_+)\times L^1(\R_+^2),$$
and the subspaces $Y_{\tau}, Y_{\eta}\subset L^1(\R_+^2)$ by
$$Y_k=\left\{\varphi\in L^1(\R_+^2): \frac{\partial \varphi}{\partial a}\in L^1(\R_+), \frac{\partial \varphi}{\partial k} \in L^1(\R_+)\right\},$$
for $k\in\{\tau,\eta\}$. Note that if we define the norms $\|\cdot\|_{k}$ on $Y_k$ by
$$\|\varphi\|_{k}=\|\varphi\|_{L^1(\R_+^2)}+\|\partial_a \varphi\|_{L^1(\R_+^2)}+\|\partial_{k}\varphi\|_{L^1(\R_+^2)}; \ \forall \varphi\in Y_k,$$ 
then $(Y_k,\|\cdot\|_k)$ becomes a Banach space. Denoting by $\Co^1_{c}(\R_+^2)$ the set of all compactly supported $\mathcal{C}^1$-function in $\R_+^2$, we see that $\Co^1_{c}(\R_+^2)\subset Y_k$ for each $k\in\{\tau,\eta\}$. Since $\Co^1_{c}(\R_+^2)$ is dense in $L^1(\R_+^2)$, then $Y_{\tau}$ and $Y_{\eta}$ are also dense in $L^1(\R_+^2)$. Moreover, for each $\varphi\in \Co^1_{c}(\R_+^2)\cap Y_k$, we have
$$\|\varphi(0,\cdot)\|_{L^1(\R_+)}\leq \|\partial_a \varphi\|_{L^1(\R_+^2)}, \quad \textnormal{ and } \quad \|\varphi(\cdot,0)\|_{L^1(\R_+)}\leq \|\partial_{k}\varphi\|_{L^1(\R_+^2)},$$
for $k\in\{\tau,\eta\}$. With these estimates, the following trace lemma holds true.

\begin{lemma}
There exists a unique linear operator $\Pi_{i,k}\in \L(Y_{k},L^1(\R_+))$ for each $i=1,2$ and $k\in\{\tau,\eta\}$ such that for all $\varphi_k\in \Co^1_{c}(\R_+^2)\cap Y_k:$
$$\Pi_{1,k} \varphi_k=\varphi_k(0,\cdot), \qquad \Pi_{2,k} \varphi_k=\varphi_k(\cdot,0).$$
\end{lemma}
We can now define the linear operators $A_{h,2}:D(A_{h,2})\subset X_2\to X_2$, $A_{h,3}:D(A_{h,3})\subset X_2\to X_2$ and $A_{m,2}:D(A_{m,2})\subset X_2\to X_2$ by
$$A_{h,2}\begin{pmatrix}
0_{L^1(\R_+)} \\
\phi
\end{pmatrix}
=\begin{pmatrix}
-\Pi_{2,\tau} \phi \\
-\partial_a \phi -\partial_{\tau} \phi -(\mu_h+\nu_h+\gamma_h)\phi
\end{pmatrix}, \quad A_{h,3}\begin{pmatrix}
0_{L^1(\R_+)} \\
\phi
\end{pmatrix}
=\begin{pmatrix}
-\Pi_{2,\eta} \phi \\
-\partial_a \phi -\partial_{\eta} \phi -(\mu_h+k_h)\phi
\end{pmatrix},$$
and
$$A_{m,2}\begin{pmatrix}
0_{L^1(\R_+)} \\
\phi
\end{pmatrix}
=\begin{pmatrix}
-\Pi_{2,\tau} \phi \\
-\partial_a \phi -\partial_{\tau} \phi -(\mu_m+\nu_m)\phi
\end{pmatrix},$$
where $D(A_{h,2})=D(A_{m,2})=\{0_{L^1(\R_+)}\}\times (Y_\tau\cap \ker(\Pi_{1,\tau}))$ and $D(A_{h,3})=\{0_{L^1(\R_+)}\}\times (Y_\eta\cap \ker(\Pi_{1,\eta}))$. 

We find that for $\lambda\in \C$ such that $\Re(\lambda)\geq -\mu_0$, we have $\lambda\in \rho(A_{h,2})\cap\rho(A_{h,3})\cap \rho(A_{m,2})$, and for each $(\psi_1,\psi_2)^T\in X_2$ we have the following explicit formula for the resolvent:
\begin{equation}
\left\{
\begin{array}{ll}
\label{Eq:Resolv2}
(\lambda I-A_{h,2})^{-1}\begin{pmatrix}
\psi_1 \\
\psi_2
\end{pmatrix}
=\begin{pmatrix}
0_{L^1(\R_+)}\\
\phi_1
\end{pmatrix}, \
(\lambda I-A_{h,3})^{-1}\begin{pmatrix}
\psi_1 \\
\psi_2
\end{pmatrix}
=\begin{pmatrix}
0_{L^1(\R_+)}\\
\phi_2
\end{pmatrix}, \\
(\lambda I-A_{m,2})^{-1}\begin{pmatrix}
\psi_1 \\
\psi_2
\end{pmatrix}
=\begin{pmatrix}
0_{L^1(\R_+)}\\
\phi_3
\end{pmatrix}
\end{array}
\right.
\end{equation}
if and only if
\begin{flalign*}
   \phi_1(a,\tau)&= \displaystyle \mathbf{1}_{[0,a]}(\tau)\left[
\psi_1(a-\tau)\exp\left(-\int_0^\tau (\mu_h(a+s-\tau)+\nu_h(a+s-\tau,s)+\gamma_h(a+s-\tau,s))\d s\right)\right. \\
+ & \left.\displaystyle\int_0^\tau \psi_2(a+s-\tau,s)\exp\left(-\int_s^\tau (\mu_h(a+\xi-\tau,\xi)+\nu_h(a+\xi-\tau,\xi)+\gamma_h(a+\xi-\tau,\xi))\d \xi\right)\d s\right],
    \end{flalign*}
\begin{equation*}
    \begin{array}{lcl}
   \phi_2(a,\eta)&=& \displaystyle \mathbf{1}_{[0,a]}(\eta)\left[]
\psi_1(a-\eta)\exp\left(-\int_0^\eta (\mu_h(a+s-\eta)+k_h(a+s-\eta,s))\d s\right)\right. \\
&&\left.+\displaystyle \int_0^\eta \psi_2(a+s-\eta,s)\exp\left(-\int_s^\eta (\mu_h(a+\xi-\eta,\xi)+k_h(a+\xi-\eta,\xi))\d \xi\right)\d s\right],
    \end{array}
\end{equation*}
and
\begin{equation*}
    \begin{array}{lcl}
  \phi_3(a,\tau)&=& \displaystyle \mathbf{1}_{[0,a]}(\tau)\left[
\psi_1(a-\tau)\exp\left(-\int_0^\tau (\mu_m(a+s-\tau)+\nu_m(a+s-\tau,s))\d s\right)\right. \\
&&\left.+\displaystyle \int_0^\tau \psi_2(a+s-\tau,s)\exp\left(-\int_s^\tau (\mu_m(a+\xi-\tau,\xi)+\nu_m(a+\xi-\tau,\xi))\d \xi\right)\d s\right].
    \end{array}
\end{equation*}
We now consider the Banach space
$$X=X_1\times X_2\times X_2\times X_1\times X_2$$
and its positive cone
$$X_+=X_{1,+}\times X_{2,+}\times X_{2,+}\times X_{1,+}\times X_{2,+}$$
endowed with the usual product norm and define the linear operator $A:D(A)\subset X\to X$ by:
$$D(A)=D(A_{h,1})\times D(A_{h,2})\times D(A_{h,3})\times D(A_{m,1})\times D(A_{m,2}), \quad A=\diag(A_{h,1}, A_{h,2}, A_{h,3}, A_{m,1}, A_{m,2}).$$
We set $X_0=\overline{D(A)}$, the closure of $D(A)$, which is given by
$$X_0 = \left(\{0_\R\}\times L^1(\R_+)\right)\times  \left(\{0_{L^1(\R_+)}\}\times L^1(\R_+^2)\right)^2\times\left(\{0_\R\}\times L^1(\R_+)\right)\times\left( \{0_{L^1(\R_+)}\}\times L^1(\R_+^2)\right),$$
so that $D(A)$ is not dense in $X$, and we define its positive cone $X_{0,+}=X_0\cap X_+$. Importantly, note that the nonlinear part of \eqref{Eq:Model} is not well defined on $X_0$ due to the term $N_h(t)$ and consequently is not locally Lipschitz continuous. To fix this, we set 
\begin{equation}\label{Eq:u}
u(t)=\left(0_\R,
S_h(t,\cdot), 0_{L^1}, I_h(t,\cdot,\cdot), 0_{L^1},
R_h(t,\cdot,\cdot), 0_\R, S_m(t,\cdot), 0_{L^1}, I_m(t,\cdot,\cdot)\right)^T,
\end{equation}
and we define the space
$$X_{\ep}=\{u(t)\in X_0: \Tau(u(t))\geq \ep\}\subset X_0,$$
where $\Tau:X\to X$ is the operator defined by
$$\Tau(u(t))=\|S_h(t,\cdot)\|_{L^1(\R_+)}+\|I_h(t,\cdot,\cdot)\|_{L^1(\R_+^2)}+\|R_h(t,\cdot,\cdot)\|_{L^1(\R_+^2)},$$
for each $u(t)\in X$ defined by \eqref{Eq:u}. Note that $X_{\ep=0}$ corresponds to the space $X_0=\overline{D(A)}$. We can now define the nonlinear operator $F_\ep:X_\ep\to X$ by
$$F_\ep(u(t))=\begin{pmatrix}
\Lambda_h \\
-\frac{S_h(t,\cdot)}{N_h(t)}\int_0^\infty \int_0^\infty \theta \beta_m(s,\tau)I_m(t,s,\tau)\d s~\d \tau+\int_0^\infty k_h(\cdot,\eta) R_h(t,\cdot,\eta)\d \eta \\
\frac{S_h(t,\cdot)}{N_h(t)}\int_0^\infty \int_0^\infty \theta \beta_m(s,\tau)I_m(t,s,\tau)\d s~\d \tau \\
0 \\
\int_0^\infty \gamma_h(\cdot,\tau) I_h(t,\cdot,\tau)\d \tau \\
0 \\
\Lambda_m \\
-\frac{S_m(t,\cdot)}{N_h(t)}\int_0^\infty \int_0^\infty \theta  \beta_h(s,\tau)I_h(t,s,\tau)\d s~\d \tau \\
\frac{S_m(t,\cdot)}{N_h(t)}\int_0^\infty \int_0^\infty \theta  \beta_h(s,\tau)I_h(t,s,\tau)\d s~\d \tau \\
0
\end{pmatrix},$$
that is well defined for each $\ep>0$. Then \eqref{Eq:Model} rewrites as the following non-densely defined abstract Cauchy problem:
\begin{equation}\label{Eq:CauchyPb}
\frac{\d u}{\d t}(t)=A u(t)+F_\ep(u(t)), \ t>0, \ u(0)=\hat{u}_0\in X_{0,+},
\end{equation}
where
\begin{equation}\label{Eq:u_hat}
\hat{u}_0=\left(0,
S_{h,0}, 0, I_{h,0}, 0,
R_{h,0}, 0, S_{m,0}, 0, I_{m,0}\right)^T,
\end{equation}
for $\ep>0$. 

\subsection{Well-posedness}

Using the above semigroup formulation, we now give the proof of Theorem \ref{Thm:existence}. The proof is split into several steps.

\textbf{Step 1}: we start by proving the existence of positive solutions on some interval $[0,\tau)$. From \eqref{Eq:Resolv1}-\eqref{Eq:Resolv2}, we see that  $(A,D(A))$ is a closed linear operator such that $(-\mu_0,\infty)\subset \rho(A)$ with $\overline{D(A)}=X_0$.  Moreover, using the explicit formula \eqref{Eq:Resolv1}-\eqref{Eq:Resolv2} of the resolvent, we obtain
$$\|(\lambda I-A_{k,1})^{-1}\|_{X_1}\leq \frac{1}{(\lambda+\mu_0)}, \quad \|(\lambda I-A_{k,2})^{-1}\|_{X_2}\leq \frac{1}{(\lambda+\mu_0)}, \quad  \|(\lambda I-A_{h,3})^{-1}\|_{X_2}\leq \frac{1}{(\lambda+\mu_0)}$$
for each $k\in\{h,m\}$. It then follows that $A$ is a Hille-Yosida operator and generates a locally Lipschitz continuous integrated semigroup, denoted by $\{S_A(t)\}_{t\geq 0}\subset \L(X)$. We also clearly see that $A$ is resolvent positive, \textit{i.e.} $S_A(t) X_+\subset X_+$, whence the semigroup $\{S_A(t)\}_{t\geq 0}$ is positive. The rest of the proof consists on a fixed point argument. First, let $u_0=(S_{h,0}, I_{h,0}, R_{h,0}, S_{m,0}, I_{m,0})\in \X_+\cap \X_{\overline{\ep}}$, $m=2\|u_0\|_{\X}$ and $\ep=\overline{\ep}/2$. Then consider the constant $\omega \geq 0$ defined by
$$\omega=\frac{m\theta }{\ep}\max\{\|\beta_h\|_{L^\infty},\|\beta_m\|_{L^\infty}\}.$$
It follows that the linear operator $A^\omega$ defined by $A^\omega=A-\omega I$ is also a Hill-Yosida operator and generates a locally Lipschitz continuous integrated semigroup, denoted by $\{S_{A^\omega}(t)\}_{t\geq 0}\subset \L(X)$ which is positive. Now define the set
$$B_m=\{u\in X: \|u\|_{X}\leq m\}.$$
We can readily check that $F_{\ep}\in \Co^\infty(X_{\ep},X)$ and there exists $k>0$ such that for every $(u, u_1, u_2)\in (X_{\ep}\cap X_+\cap B_m)^3$, we have the following Lipschitz and positivity properties
$$\|F_{\ep}(u_1)-F_{\ep}(u_2)\|_{X}\leq k\|u_1-u_2\|_{X}, \qquad F^\omega(u):=F_{\ep}(u)+\omega u\in X_+.$$
Now, let the constant
$$\tau=\min\left\{\frac{1}{2(k+\omega)},\frac{\ln(2)}{\|\mu_h\|_{L^\infty}+\|\nu_h\|_{L^\infty}+\|k_h\|_{L^\infty}+\|\gamma_h\|_{L^\infty}+\omega}\right\}>0,$$
and define $\{T_{(A^\omega)_0}(t)\}_{t\geq 0}$ the $\Co_0$-semigroup generated by the linear operator $(A^\omega)_0:D(A_0)\subset X\to X$, that is the part of $A^\omega$ in $X_0$. It follows that $\|T_{(A^\omega)_0}(t)u\|_{X}\leq \|u\|_{X}e^{-(\mu_0+\lambda_{m\to h})t}$ for each $t\geq 0$ and $u\in X$. Let the space
$$\Zz :=\Co^0([0,\tau],X_{\ep}\cap X_+\cap B_m),$$
be equipped with the metric
$$d(u_1,u_2)=\max_{t\in[0,\tau]}\left(|u_1(t)-u_2(t)|\right), \quad \forall (u_1,u_2)\in \Zz ^2.$$
Let $G:\Zz \to \Co^0\left([0,\tau],X)\right)$ be the operator defined by
$$G(u)(t)={T_{(A^\omega)_0}(t)}\hat{u}_0+\frac{\d}{\d t}\left(S_{A^\omega}* F^\omega(u)\right)(t),$$
with $\hat{u}_0\in X_{\overline{\ep}}$ defined by \eqref{Eq:u_hat} and where $*$ denotes the convolution product, \textit{i.e.}
$$(S_{A^\omega}* F^\omega(u))(t)=\int_0^t S_{A^\omega}(t-s)F^\omega(u(s))\d s.$$
Since $F^\omega(u)\in L^1((0,\tau),X)$, it follows by the Kellermann-Hieber theorem \cite{KellermanHieber89} (see also \cite[Theorem 3.2, p. 133]{MagalRuan2018}) that the map $t\longmapsto (S_{A^\omega}* F^\omega(u))(t)$ is continuously differentiable and satisfies:
$$\left\|\frac{\d}{\d t}(S_{A^\omega}* F^\omega(u))(t)\right\|_{X}\leq \int_0^t e^{-(\mu_0+\omega)(t-s)}\|F^\omega(u)(s)\|_{X} \d s\leq \tau (k+\omega)\max_{s\in[0,\tau]}\|u(s)\|_{X}\leq \tau m(k+\omega),$$
for each $t\in[0,\tau]$. The fact that $\{T_{(A^\omega)_0}(t)\}_{t\geq 0}$ is a $\Co_0$-semigroup induces that $G(\Zz)\subset \Co([0,\tau],X)$ and $G$ is well-defined. Moveover, by definitions of $m$ and $\tau$ we deduce that $G(\Zz)\subset \Co([0,\tau],B_m)$. Now, using the following approximation formula \cite[Proposition 3.4.8, p. 122]{MagalRuan2018}
\begin{equation}\label{Eq:approx}
\frac{\d}{\d t}(S_{A^\omega}* F^\omega(u))(t)=\lim_{\lambda\to\infty}\int_0^t T_{(A^\omega)_0}(t-s)\lambda (\lambda I-A^\omega)^{-1}F^\omega(u(s))\d s
\end{equation}
and the positivity properties of $F^\omega$ and $\{{S_{A^\omega}(t)}\}_{t\geq 0}$, we deduce that $G(\Zz)\subset \Co([0,\tau],X_+\cap B_m)$. We now prove that $G(u)(t)\in X_{\ep}$ for each $u\in \Zz$ and $t\in[0,\tau]$. Letting the constant
$$C_h=\|\mu_h\|_{L^\infty}+\|\nu_h\|_{L^\infty}+\|k_h\|_{L^\infty}+\|\gamma_h\|_{L^\infty}+\omega.$$
we see that
$$T_{(A^\omega)_0}(t)\hat{u}_0(a,\tau)\geq \begin{pmatrix}
0 \\
\mathbf{1}_{[t,\infty)}(a)S_{h,0}(a-t)e^{-t C_h} \\
0 \\
\mathbf{1}_{[t,\infty)}(\tau)\mathbf{1}_{[0,a]}(\tau)I_{h,0}(a-t,\tau-t)e^{-t C_h} \\
0 \\
\mathbf{1}_{[t,\infty)}(\tau)\mathbf{1}_{[0,a]}(\tau)R_{h,0}(a-t,\tau-t)e^{-t C_h} \\
0 \\
0 \\
0 \\
0
\end{pmatrix}.$$
By definition of $G$, we can compute $\|\Tau(G(u)(t))$ as follows:
\begin{flalign*}
\Tau(G(u(t))\geq&\left(\int_t^\infty S_{h,0}(a-t) \d a+\int_t^\infty \int_\tau^\infty (I_{h,0}(a-t,\tau-t)+R_{h,0}(a-t,\tau-t))\d a\d \tau\right)e^{-t C_h} \\
\geq& \left(\|S_{h,0}\|_{L^1(\R_+)}+\|I_{h,0}\|_{L^1(\R_+^2)}+\|R_{h,0}\|_{L^1(\R_+^2)}\right)e^{-t C_h}\\
\geq&\overline{\ep}e^{-\tau C_h},
\end{flalign*}
for each $t\in[0,\tau]$, since $u_0\in \X_{\overline{\ep}}$. By definition of $\tau$, it comes that for each $t\in [0,\tau]$ and $u\in \Zz$ we have
$$\Tau(G(u(t)))\geq \ep.$$
Consequently we have proved that $G(u)(t)\in X_{\ep}$ for each $t\in[0,\tau]$, whence $G(\Zz)\subset \Zz$ and $G$ preserves the space $\Zz$. For each $(u_1,u_2)\in \Zz^2$, the following computations:
\begin{flalign*}
\|G(u_1)-G(u_2)\|_{\Zz }&=\max_{t\in[0,\tau]}\|G(u_1(t))-G(u_2(t))\|_{X}\\
&=\max_{t\in[0,\tau]}\left\|\frac{\d}{\d t}\left(S_{A^\omega}* (F^\omega(u_1)-F^\omega(u_2))\right)(t)\right\|_{X}\\
&\leq \tau(k+\omega)\max_{t\in[0,\tau]}\|u_1(t)-u_2(t)\|_{X} \\
&\leq \tau(k+\omega)\|u_1-u_2\|_{\Zz}\\
&\leq\frac{1}{2}\|u_1-u_2\|_{\Zz},
\end{flalign*}
induce that $G$ is a $1/2$-shrinking operator. The Banach-Picard theorem then implies the existence and uniqueness of a mild solution $u\in \Co([0,\tau), X_{\ep}\cap X_+)$ for the Cauchy problem \eqref{Eq:CauchyPb}. We remind that a mild solution of \eqref{Eq:Model} or \eqref{Eq:CauchyPb} is a continuous function $u\in\Co([0,\tau),X_0)$ such that
$$\int_0^t u(s)\d s\in D(A), \ \forall t\geq 0 \textnormal{ and } u(t)=\hat{u}_0+A\int_0^t u(s)\d s+\int_0^s F_{\ep}(u(s))\d s, \ \forall t\geq 0.$$
This solution is defined by \eqref{Eq:u} and satisfies the Volterra integral formulation as stated in Theorem \ref{Thm:existence}. 

\textbf{Step 2}: we prove Theorem \ref{Thm:existence} for initial conditions in $D(A)$. Let $\hat{u}_0\in D(A)$. The mild solution $u$ is thus clearly continuously differentiable and becomes classical: \textit{i.e.} $u\in \Co^1([0,\tau),X_{\ep}\cap X_+)$ and is solution to the PDE problem \eqref{Eq:Model}-\eqref{Eq:Bound_cond}-\eqref{Eq:Ini_cond}. After integration, we notice that the following inequality holds true
\begin{equation}\label{Eq:Deriv_Nh}
N_h'(t)=\Lambda_h-\int_0^\infty \mu_h(a)\left(S_h(t,a)+R_h(t,a)\right)\d a-\int_0^\infty \int_0^\infty \left(\mu_h(a)+\nu_h(a,\tau)\right)I_h(t,a,\tau)\d a~\d \tau.
\end{equation}
By the positivity of solutions, it follows that
$$N_h'(t)\leq \Lambda_h-\mu_0 N_h(t), \qquad N_h'(t)\geq \Lambda_h-(\|\mu_h\|_{L^\infty}+\|\nu_h\|_{L^\infty})N_h(t)$$
and inequalities \eqref{Eq:Nh_gronw1}-\eqref{Eq:Nh_gronw2} comes from the use of Gronwall inequality. It also follows that
$$\Tau(u(t))=N_h(t)\geq \Tau(\hat{u}_0)e^{-(\|\mu_h\|_{L^\infty}+\|\nu_h\|_{L^\infty})t}+\left(\frac{\Lambda_h}{\|\mu_h\|_{L^\infty}+\|\nu_h\|_{L^\infty}}\right)\left(1-e^{-(\|\mu_h\|_{L^\infty}+\|\nu_h\|_{L^\infty})t}\right)$$
for each $t\in[0,\tau)$. By assumption on $\overline{\ep}$, we know that
$$\overline{\ep}\leq \frac{\Lambda_h}{\|\mu_h\|_{L^\infty}+\|\nu_h\|_{L^\infty}+\|\gamma_h\|_{L^\infty}+\|k_h\|_{L^\infty}}\leq \frac{\Lambda_h}{\|\mu_h\|_{L^\infty}+\|\nu_h\|_{L^\infty}}$$
so there exists a constant $c\geq 0$ such that
$$\frac{\Lambda_h}{\|\mu_h\|_{L^\infty}+\|\nu_h\|_{L^\infty}}=\overline{\ep}+c$$
whence
\begin{flalign*}
\Tau(u(t))&\geq \frac{\Lambda_h}{\|\mu_h\|_{L^\infty}+\|\nu_h\|_{L^\infty}}+\left(\Tau(\hat{u}_0)-\left(\frac{\Lambda_h}{\|\mu_h\|_{L^\infty}+\|\nu_h\|_{L^\infty}}\right)\right)e^{-(\|\mu_h\|_{L^\infty}+\|\nu_h\|_{L^\infty})t} \\
&\geq \overline{\ep}+c\left(1-e^{-(\|\mu_h\|_{L^\infty}+\|\nu_h\|_{L^\infty})t}\right)+(\Tau(\hat{u}_0)-\overline{\ep})e^{-(\|\mu_h\|_{L^\infty}+\|\nu_h\|_{L^\infty})t}\geq \overline{\ep}
\end{flalign*}
since $\Tau(\hat{u}_0)\geq \overline{\ep}$ by assumption. We then deduce that $u\in \Co^1([0,\tau), X_{\overline{\ep}}\cap X_+)$. From here, we define the operator $\tilde{G}:\tilde{Z}\to \Co^0([0,\tau],X))$ where $$\tilde{G}(u)(t)={T_{(A^\omega)_0}(t)}\hat{u}_0+\frac{\d}{\d t}\left(S_{A^\omega}* F^\omega(u)\right)(t),$$
for each $u\in \tilde{Z}:=\Co^0([0,\tau],X_{\overline{\ep}}\cap X_+\cap B_m)$. Proceeding as in Step 1, we can show that $G$ is a $1/2$-shrinking operator with $G(\tilde{Z})\subset \tilde{Z}$. Since $u_0\in X_{\overline{\ep}}$, we can readily use some standard time extending properties of the solution to extend the solution $u$ over a maximal interval $[0,t_{\max})$ with $t_{\max}>0$. We remark that, using integrations, the inequality
\begin{equation}\label{Eq:Deriv_Nv}
N_m'(t)=\Lambda_m-\int_0^\infty \mu_m(a)S_m(t,a)\d a-\int_0^\infty \int_0^\infty \left(\mu_m(a)+\nu_m(a,\tau)\right)I_m(t,a,\tau)\d a~\d \tau,
\end{equation}
is satisfied, which induces the fact that
$$N_m'(t)\leq \Lambda_m-\mu_0 N_m(t), \qquad N_m'(t)\geq \Lambda_m-\left(\|\mu_m\|_{L^\infty}+\|\nu_m\|_{L^\infty}\right)N_m(t),$$
whence the inequalities \eqref{Eq:Nv_gronw1}-\eqref{Eq:Nv_gronw2} hold true. We deduce that
\begin{equation}\label{Eq:ineq:Nh-Nv}
\limsup_{t\to \infty}N_h(t)\leq \frac{\Lambda_h}{\mu_0}, \qquad \limsup_{t\to\infty}N_m(t)\leq \frac{\Lambda_m}{\mu_0},
\end{equation}
and that the inequalities \eqref{Eq:bounded_human}-\eqref{Eq:bounded_mector} are satisfied by positivity of the solutions. Using these estimates, it readily follows (see \textit{e.g.} \cite[Theorem 6.1.4, p. 185]{Pazy83}) that the solution is global, \textit{i.e.} $t_{\max}=\infty$ and $u\in\Co^1(\R_+, X_{\overline{\ep}}\cap X_+)$. Finally, from \eqref{Eq:Model} we see that
$$\left(\frac{\partial}{\partial t}+\frac{\partial}{\partial a}\right)S_m(t,a)\geq-\left(\mu_m(a) +\theta \|\beta_h\|_{L^\infty(\R_+^2)}\right)S_m(t,a), \quad S_m(t,0)=\Lambda_m, \quad \forall (t,a)\in \R_+^2$$
which leads to \eqref{Eq:bounded_Sv}. Moreover, we see on one hand that
$$\left(\frac{\partial}{\partial t}+\frac{\partial}{\partial a}\right)S_h(t,a)\geq-\left(\mu_h(a) +\frac{N_m(t)}{N_h(t)}\theta \|\beta_m\|_{L^\infty(\R_+^2)}\right)S_h(t,a),  \quad S_h(t,0)=\Lambda_h, \quad \forall (t,a)\in \R_+^2.$$
On the other hand, we deduce from \eqref{Eq:Nh_gronw2}-\eqref{Eq:Nv_gronw1} that
$$\limsup_{t\to \infty}N_h(t)\geq \frac{\Lambda_h}{\|\mu_h\|_{L^\infty}+\|\nu_h\|_{L^\infty}}, \qquad \liminf_{t\to\infty}N_m(t)\leq \frac{\Lambda_m}{\mu_0}.$$
These two latter points combined together imply that \eqref{Eq:bounded_Sh} is satisfied.

\textbf{Step 3}: we now prove Theorem \ref{Thm:existence}. Suppose now that $\hat{u}_0\in X_0$ so that $u\in \Co([0,\tau),X_{\ep}\cap X_+)$ is a mild solution to \eqref{Eq:CauchyPb}. Since $\overline{D(A)}=X_0$, it follows that there exists a sequence of initial conditions $\{\hat{u}^k_0\}_{k\geq 0}\subset D(A)^\N$ such that $\lim_{k\to \infty}\|(\hat{u}^k_0-\hat{u}_0)\|_{X}=0$. For each $k\geq 0$, there exists a unique solution $u_k\in \Co^1(\R_+,X_{\overline{\ep}}\cap X_+)$ to \eqref{Eq:Model}-\eqref{Eq:Bound_cond}-\eqref{Eq:Ini_cond} with initial condition $\hat{u}^k_0$. For each $t\in[0,\tau)$ we can compute
\begin{flalign*}
\|u(t)-u^k(t)\|_{X}&=\left\|S_{(A^\omega)_0}(t)\hat{u}_0-S_{(A^\omega)_0}(t)\hat{u}^k_0+\frac{\d}{\d t}(S_{A^\omega}*F^\omega(u))(t)-\frac{\d}{\d t}(S_{A^\omega}*F^\omega(u^k))(t)\right\|_{X} \\
& \leq \|\hat{u}_0-\hat{u}^k_0\|_{X}+\frac{1}{2}\max_{s\in[0,\tau]}\|u(s)-u^k(s)\|_{X},
\end{flalign*}
whence
$$\|u-u^k\|_{\Zz}=\max_{t\in [0,\tau]}\|u(t)-u^k(t)\|_{X}\leq 2\|\hat{u}_0-\hat{u}^k_0\|_{X}\underset{k\to 0}{\to}0.$$
Writing
$$u(t)=u^k(t)+u(t)-u^k(t),$$
we see that
$$\Tau(u(t))\geq \Tau(u^k(t))-\|u(t)-u^k(t)\|_{X}.$$
On one hand we obtain
$$\Tau(u(t))\geq \overline{\ep}-\|u(t)-u^k(t)\|_{X},$$
and on the other hand we have
\begin{flalign*}
\Tau(u(t))\geq& \Tau(\hat{u}^k_0)e^{-(\|\mu_h\|_{L^\infty}+\|\nu_h\|_{L^\infty})t}+\frac{\Lambda_h}{\|\mu_h\|_{L^\infty}+\|\nu_h\|_{L^\infty}}\left(1-e^{-(\|\mu_h\|_{L^\infty}+\|\nu_h\|_{L^\infty})t}\right)-\|u(t)-u^k(t)\|_{X}\\
\geq& \Tau(\hat{u}_0)e^{-(\|\mu_h\|_{L^\infty}+\|\nu_h\|_{L^\infty})t}+\frac{\Lambda_h}{\|\mu_h\|_{L^\infty}+\|\nu_h\|_{L^\infty}}\left(1-e^{-(\|\mu_h\|_{L^\infty}+\|\nu_h\|_{L^\infty})t}\right)-\|u(t)-u^k(t)\|_{X}. \\
&-\left\|\Tau(\hat{u}_0-\hat{u}^k_0)\right\|_{X}e^{-(\|\mu_h\|_{L^\infty}+\|\nu_h\|_{L^\infty})t}.
\end{flalign*}
Letting $k$ goes to infinity, it respectively follows that $\Tau(u(t))\geq \overline{\ep}$ whence $u\in \Co([0,\tau), X_{\overline{\ep}}\cap X_+$, and the inequality \eqref{Eq:Nh_gronw1} holds for each $t\in[0,\tau)$ since we have
$$\|\Tau(\hat{u}_0-\hat{u}^k_0)\|_{X}\leq \|\hat{u}_0-\hat{u}^k_0\|_{X}\underset{k\to\infty}{\to}0.$$
Similarly, from
\begin{flalign*}
\Tau(u(t))\leq& \Tau(u^k(t))+\|u(t)-u^k(t)\|_{X}\\
\leq& \Tau(\hat{u}^k_0)e^{-\mu_0 t}+\frac{\Lambda_h}{\mu_0}\left(1-e^{-\mu_0 t}\right)+\|u(t)-u^k(t)\|_{X}\\
\leq& \Tau(\hat{u}_0)e^{-\mu_0 t}+\Tau(\hat{u}_0-\hat{u}^k_0)e^{-\mu_0 t}+\frac{\Lambda_h}{\mu_0}\left(1-e^{-\mu_0 t}\right)+\|u(t)-u^k(t)\|_{X},
\end{flalign*}
we see that \eqref{Eq:Nh_gronw2} holds for each $t\in[0,\tau)$. The inequalities \eqref{Eq:Nv_gronw1}-\eqref{Eq:Nv_gronw2} are proved similarly by using the operator $\overline{\Tau}:X\to X$ defined by
$$\overline{\Tau}(u(t))=\|S_m(t,\cdot)\|_{L^1(\R_+)}+\|I_m(t,\cdot,\cdot)\|_{L^1(\R_+^2)},$$
for each $u(t)\in X$ defined by \eqref{Eq:u}. 
As in Step 2, we deduce that the solution $u$ is global, \textit{i.e.} $u\in \Co(\R_+,X_{\overline{\ep}}\cap X_+)$ and the above estimates then hold for each $t\geq 0$. It follows that $u$ satisfies \eqref{Eq:ineq:Nh-Nv} and \eqref{Eq:bounded_human}-\eqref{Eq:bounded_mector} hold true by positivity. Finally, from these estimates and using the Volterra integral formulation of the solution $u$ as stated in Theorem \ref{Thm:existence}, we can show that the inequalities \eqref{Eq:bounded_Sh}-\eqref{Eq:bounded_Sv} hold, which ends the proof of Theorem \ref{Thm:existence}.

\section{The basic reproduction number and proof of Theorem \ref{Thm:stab}}
Here we derive the basic reproduction number $\mathcal{R}_0$ of Model \eqref{Eq:Model} and give details on the stability results of the disease-free equimibrium $E^0$.

Let $\ep\geq 0$. The linearised system of \eqref{Eq:Model} around $E_0$ is:
\begin{equation}\label{Eq:Linearized}
    \left\{
    \begin{array}{rcl}
        \frac{\d u}{\d t}(t)&=&A u(t)+D_{E_0}F_{\ep}(u(t))\\
         u(0)&=&\hat{u}_0\in \overline{D(A)}
    \end{array}
    \right.
\end{equation}
where $D_{E_0}F_{\ep}:X_{0}\to X$ denotes the differential of $F_{\ep}$ around $E_0$ and is defined by
\begin{equation*}
  D_{E_0}F_{\ep}\begin{pmatrix}
    0 \\
    S_h \\
    0 \\
    I_h \\
    0 \\
    R_h \\
    0 \\
    S_m \\
    0 \\
    I_m
    \end{pmatrix}=
    \begin{pmatrix}
    0 \\
    -\frac{S_h^0(a)}{N_h^0}\int_0^\infty \int_0^\infty \theta \beta_m(s,\tau)I_m(s,\tau)\d s~\d \tau \vspace{0.1cm} \\
  \frac{S_h^0(a)}{N_h^0}\int_0^\infty \int_0^\infty \theta \beta_m(s,\tau)I_m(s,\tau)\d s~\d \tau\\
    0 \\
   0 \\
    0 \\
    0 \\
    - \frac{S_m^0(a)}{N_h^0}\int_0^\infty \int_0^\infty \theta \beta_h(s,\tau)I_h(s,\tau)\d s~\d \tau \\
  \frac{S_m^0(a)}{N_h^0}\int_0^\infty \int_0^\infty \theta \beta_h(s,\tau)I_h(s,\tau)\d s~\d \tau \\
    0
    \end{pmatrix}+
    \begin{pmatrix}
    0 \\
   \int_0^\infty k_h(a,\eta)R_h(a,\eta) \vspace{0.1cm} \\
  0 \\
    0 \\
   \int_0^\infty \gamma_h(a,\tau) I_h(a,\tau)\d \tau \\
    0 \\
    0 \\
   0 \\
  0 \\
    0
    \end{pmatrix}
\end{equation*}
that we decompose into
$$D_{E_0}F_{\ep}=(D_{E_0}F_{\ep})_1+(D_{E_0}F_{\ep})_2$$
with
$$N_h^0=\int_0^\infty S_h^0(a)\d a=\Lambda_h \int_0^\infty e^{-\int_0^a \mu_h(s)\d s}\d a.$$
Note that since $N_h^0>0$, then the differential operator $D_{E_0}F_{\ep}$ is well-defined on $X_0=\overline{D(A)}$ whatever $\ep\geq 0$. 

\subsection{The basic reproduction number $\mathcal{R}_0$ of Model \eqref{Eq:Model}}

We derive the $\mathcal{R}_0$ in several steps. 

\textbf{Step 1:} we begin by computing the next generation operator (see \cite{Diekmann90, Inaba2012}). Let
$$(0,S_h,0,I_h,0,R_h,0,S_v,0,I_v)^T\in X_0$$
be a solution of \eqref{Eq:Linearized}, that is the linearized system of the model \eqref{Eq:Model} around the disease-free equilibrium $E_0$. It follows that $(I_h, I_m)$ satisfy the following equations:
\begin{equation*}
\left\{
\begin{array}{l}
\left(\frac{\partial}{\partial t}+\frac{\partial}{\partial a}+\frac{\partial}{\partial \tau}\right)I_h(t,a,\tau)=-\left(\mu_h(a)+\nu_h(a,\tau)+\gamma_h(a,\tau)\right)I_h(t,a,\tau), \vspace{0.1cm} \\
\left(\frac{\partial }{\partial t}+\frac{\partial }{\partial a}+\frac{\partial }{\partial \tau}\right)I_m(t,a,\tau)=-(\mu_m(a)+\nu_m(a,\tau))I_m(t,a,\tau),\\
I_h(t,a,0)=B_h(t,a), \qquad I_h(t,0,\tau)=0, \qquad I_h(0,a,\tau)=I_{h,0}(a,\tau),\\
I_m(t,a,0)=B_m(t,a), \qquad I_m(t,0,\tau)=0, \qquad I_m(0,a,\tau)=I_{m,0}(a,\tau)
\end{array}
\right.
\end{equation*}
where $B_h$ and $B_m$ respectively denote the number of newly infected humans and mosquitoes, and are defined as:
\begin{equation*}
\left\{
\begin{array}{rcl}
B_h(t,a)&=  \frac{\pi_h(a)}{\int_0^\infty \pi_h(s)\d s} \int_0^\infty \int_0^\infty \theta \beta_m(s,\tau) I_m(t,s,\tau) \d s \d \tau , \vspace{0.1cm}\\ 
B_m(t,a)&=  \frac{\Lambda_m \pi_m(a)}{\Lambda_h \int_0^\infty \pi_h(s)\d s} \int_0^\infty \int_0^\infty \theta \beta_h(s,\tau) I_h(t,s,\tau) \d s \d \tau,
\end{array}
\right.
\end{equation*}
with $\pi_h(a)= e^{-\int_0^a \mu_h(s)\d s}$ and $\pi_m(a)= e^{-\int_0^a \mu_m(s)\d s}$ the survival probabilities, from birth until age $a$, for humans and mosquitoes respectively, in absence of disease. As a result, $(I_h,I_m)$ are given by the following Volterra integral formulation:
\begin{equation*}
\begin{split}
I_h(t,a,\tau)= & \left\{
\begin{array}{lll} B_h(t-\tau,a-\tau)
e^{-\int_0^\tau (\mu_h(s+a-\tau)+\nu_h(s+a-\tau,s)+\gamma_h(s+a-\tau,s))\d s}, \quad \forall t>\tau, \quad \forall a\geq \tau, \\
I_{h,0}(a-t,\tau-t)e^{-\int_{\tau-t}^\tau (\mu_h(s+a-\tau)+\nu_h(s+a-\tau,s)+\gamma_h(s+a-\tau,s))\d s}, \quad \forall a\geq \tau\geq t,
\end{array}
\right.\\
I_m(t,a,\tau)=& \left\{
\begin{array}{lll}
B_m(t-\tau,a-\tau) e^{-\int_0^\tau (\mu_m(s+a-\tau)+\nu_m(s+a-\tau,s))\d s}, \quad \forall t>\tau, \quad \forall a\geq \tau, \\
I_{m,0}(a-t,\tau-t)e^{-\int_{\tau-t}^\tau (\mu_m(s+a-\tau)+\nu_m(s+a-\tau,s))\d s}, \quad \forall a\geq \tau\geq t.
\end{array}
\right.
\end{split}
\end{equation*}
Therefore, $B_h$ and $B_m$ rewrite as
\begin{flalign*}
B_h(t,a)=& \frac{\pi_h(a)}{\int_0^\infty \pi_h(s)\d s} \int_0^t \int_\tau^\infty \theta \beta_m(s,\tau)  e^{-\int_0^\tau (\mu_m(\sigma+s-\tau)+\nu_m(\sigma+s-\tau,\sigma))\d \sigma} B_m(t-\tau,s-\tau) \d s~\d \tau + f_{0,h}(t,a)\\
=&\frac{\pi_h(a)}{\int_0^\infty \pi_h(s)\d s} \int_0^t \int_0^\infty \theta \beta_m(\tau+\xi,\tau)  e^{-\int_0^\tau (\mu_m(\sigma+\xi)+\nu_m(\sigma+\xi,\sigma))\d \sigma} B_m(t-\tau,\xi) \d \xi~\d \tau+ f_{0,h}(t,a)
\end{flalign*}
and 
\begin{flalign*}
&B_m(t,a)\\
=&\frac{\Lambda_m \pi_m(a)}{\Lambda_h \int_0^\infty \pi_h(s)\d s} \int_0^t \int_\tau^\infty \theta \beta_h(s,\tau) e^{-\int_0^\tau (\mu_h(\sigma+s-\tau)+\nu_h(\sigma+s-\tau,\sigma)+\gamma_h(\sigma+s-\tau,\sigma))\d \sigma}   B_h(t-\tau,s-\tau)  \d s~\d \tau\\
&+ f_{0,m}(t,a) \\
=&\frac{\Lambda_m \pi_m(a)}{\Lambda_h \int_0^\infty \pi_h(s)\d s} \int_0^t \int_0^\infty \theta \beta_h(\tau+\xi,\tau) e^{-\int_0^\tau (\mu_h(\sigma+\xi)+\nu_h(\sigma+\xi,\sigma)+\gamma_h(\sigma+\xi,\sigma))\d \sigma}   B_h(t-\tau,\xi)  \d \xi~\d \tau+ f_{0,m}(t,a)
\end{flalign*}
where $f_{0,h}$ and $ f_{0,m}$ encounter for the initial data. Now, we see that the functions
$$(\tilde{B}_h(t), \tilde{B}_m(t))=(B_h(t,\cdot),B_m(t,\cdot))\in L^1(\R_+)^2$$
satisfy the following equations:
\begin{equation*}
\left\{
\begin{array}{rcl}
\tilde{B}_h(t)&=&\displaystyle\int_0^t G_h(\tau)(\tilde{B}_m(t-\tau))\d\tau+f_{0,h}(t,\cdot), \\
\tilde{B}_m(t)&=&\displaystyle\int_0^t G_m(\tau)(\tilde{B}_h(t-\tau))\d\tau+f_{0,m}(t,\cdot)
\end{array}
\right.
\end{equation*}
where the linear operators $G_h(\tau)$ and $ G_m(\tau)$ belong to $\L(L^1(\R_+))$, and are defined by:
\begin{flalign*}
(G_h(\tau)B)(a)&=\frac{\pi_h(a) \theta }{\int_0^\infty \pi_h(s)\d s} \int_0^\infty\beta_m(\tau+\xi,\tau)  e^{-\int_0^\tau (\mu_m(\sigma+\xi)+\nu_m(\sigma+\xi,\sigma))\d \sigma} B(\xi)\d \xi \\
(G_m(\tau)B)(a)&=\frac{\Lambda_m \pi_m(a)\theta}{\Lambda_h \int_0^\infty \pi_h(s)\d s}\int_0^\infty   \beta_h(\tau+\xi,\tau)  e^{-\int_0^\tau (\mu_h(\sigma+\xi)+\nu_h(\sigma+\xi,\sigma)+\gamma_h(\sigma+\xi,\sigma))\d \sigma} B(\xi)\d \xi
\end{flalign*}
for each $(a,\tau)\in \R_+^2$ and $B\in L^1(\R_+)$. Note that these operators $G_h(\tau)$ and $G_m(\tau)$ are called the \textit{net reproduction operators} (see \cite{Inaba2012}), that map the density of newborns to the density of their children produced at $\tau$ time later. It then follows (see \cite{Diekmann90, Inaba2012}), that the next generation operator is given by:
\begin{flalign*}
\G\begin{pmatrix}
B_1 \\
B_2
\end{pmatrix}(a)&=\displaystyle\begin{pmatrix}
\int_0^\infty (G_h(\tau)B_2)(a)\d \tau\\
\int_0^\infty (G_m(\tau)B_1)(a)\d \tau
\end{pmatrix}\\
&=\displaystyle\begin{pmatrix}
\frac{\pi_h(a)\theta }{\int_0^\infty \pi_h(s)\d s}  \int_0^\infty \int_0^\infty \beta_m(\tau+\xi,\tau)  e^{-\int_0^\tau (\mu_m(\sigma+\xi)+\nu_m(\sigma+\xi,\sigma))\d \sigma}  B_2(\xi) \d \xi~\d \tau\\
\frac{\Lambda_m \pi_m(a)\theta}{\Lambda_h \int_0^\infty \pi_h(s)\d s} \int_0^\infty \int_0^\infty  \beta_h(\tau+\xi,\tau) e^{-\int_0^\tau (\mu_h(\sigma+\xi)+\nu_h(\sigma+\xi,\sigma)+\gamma_h(\sigma+\xi,\sigma))\d \sigma}   B_1(\xi)  \d \xi~\d \tau
\end{pmatrix}.
\end{flalign*}

\textbf{Step 2:} from the next generation operator $\G$, we deduce that that the basic reproduction number $\RR_0$ is defined by the spectral radius of $\G$, denoted by $r_\sigma(\G)$. In order to compute this spectral radius, we define the operator $\G^2=\G \circ \G \in \L((L^1(\R_+))^2)$. First we know that
\begin{equation}\label{Eq:R0_Gsquare}
\RR_0=r_\sigma(\G)=\sqrt{r_\sigma(\G^2)}.
\end{equation}
Consequently, it remains to compute $r_\sigma(\G^2)$. We see that
$$\G^2(B_1,B_2)^T=(H_h(B_1),H_m(B_2))^T$$
for each $(B_1,B_2)\in L^1(\R_+)^2$, with $H_h$ and $H_m$ the linear operators on $L^1(\R_+)$ respectively by:
\begin{flalign*}
H_h(B)(a)=\frac{\Lambda_m \pi_h(a)\theta^2 }{\Lambda_h(\int_0^\infty \pi_h(s)\d s)^2} 
\left(\int_0^\infty \int_0^\infty  \beta_h(\tau+\xi,\tau) e^{-\int_0^\tau (\mu_h(\sigma+\xi)+\nu_h(\sigma+\xi,\sigma)+\gamma_h(\sigma+\xi,\sigma))\d \sigma}   B(\xi)  \d \xi~\d \tau\right)\\
\times\left(\int_0^\infty \int_0^\infty \beta_m(\tau+\xi,\tau)  e^{-\int_0^\tau (\mu_m(\sigma+\xi)+\nu_m(\sigma+\xi,\sigma))\d \sigma}  \pi_m(\xi) \d \xi~\d \tau\right) 
\end{flalign*}
and
\begin{flalign*}
H_m(B)(a)=\frac{\Lambda_m \pi_m(a)\theta^2 }{\Lambda_h(\int_0^\infty \pi_h(s)\d s)^2}  \left(\int_0^\infty \int_0^\infty  \beta_h(\tau+\xi,\tau) e^{-\int_0^\tau (\mu_h(\sigma+\xi)+\nu_h(\sigma+\xi,\sigma)+\gamma_h(\sigma+\xi,\sigma))\d \sigma}   \pi_h(\xi)  \d \xi~\d \tau\right)\\
\times \left(\int_0^\infty \int_0^\infty \beta_m(\tau+\xi,\tau)  e^{-\int_0^\tau (\mu_m(\sigma+\xi)+\nu_m(\sigma+\xi,\sigma))\d \sigma}  B(\xi) \d \xi~\d \tau\right).
\end{flalign*}
It follows from \eqref{Eq:R0_Gsquare} that
\begin{equation}\label{Eq:R0_h}
\RR_0=\sqrt{r_\sigma(\G^2)}=\sqrt{\max\{r_\sigma(H_h), r_\sigma(H_m)\}}
\end{equation}

\textbf{Step 3:} now we compute $r_\sigma(H_h)$. For this end, let us define the sets
\begin{equation}\label{Eq:omega_h}
\Omega_h=\{\xi\geq 0: \int_0^\infty \beta_h(\tau+\xi,\tau)e^{-(\|\mu_h\|_{L^\infty}+\|\nu_h\|_{L^\infty}+\|\gamma_h\|_{L^\infty})\tau}\d \tau>0\}
\end{equation}
and
\begin{equation}\label{Eq:omega_m}
\Omega_m=\{\xi\geq 0: \int_0^\infty \beta_m(\tau+\xi,\tau)e^{-(\|\mu_m\|_{L^\infty}+\|\nu_m\|_{L^\infty})\tau}\d \tau>0\}.
\end{equation}
Then, under Assumption \ref{Assum:beta}, we have $$\Omega_h \neq\emptyset, \qquad\Omega_m \neq\emptyset.$$

Next, we define the restriction of $H_h$ to $L^1(\Omega_h)$ denoted by $\tilde{H}_h\in \L(L^1(\Omega_h))$:
$$\tilde{H}_h(\tilde{B})(s)=\mathbf{1}_{\Omega_h}(s)H_h(B)(s)$$
for each $s\in \Omega_h$, $\tilde{B}\in L^1(\Omega_h)$, where $\Omega_h$ satisfies \eqref{Eq:omega_h} and with
$$B(s)=\begin{cases}
\tilde{B}(s) & \text{a.e.} \quad s\in \Omega_h,\\
0 & \text{else}.
\end{cases}$$
From Assumption \ref{Assump}, it is clear that $\tilde{H}_h$ is a compact and positive operator on $L^1(\Omega_h)$. Moreover, since $\Omega_h \neq \emptyset $, it follows that $\tilde{H}_h$ is irreducible, \textit{i.e.}
$$\tilde{H}_h(f)(s)>0 \quad \mathnormal{a.e.} \ s\in \Omega_h, \ \forall f\in L^1_+(\Omega_h)\setminus\{0\}$$
that is, sends the positive cone $L^1_+(\R_+)$ on the subset of $L^1_+(\R_+)$ of functions almost everywhere strictly positive. It follows from \cite[Theorem 3]{Pagter86} that its spectral radius is positive, that is: $r_\sigma(\tilde{H}_h)>0$. We now observe that
$$H_h(\pi_h)=\lambda_0 \pi_h$$
so that $\lambda_0$ is an eigenvalue of $H_h$ associated to the eigenfunction $\pi_h$, so that
\begin{equation}\label{Eq:spectral_rad>0}
r_\sigma(H_h)\geq \lambda_0
\end{equation}
where
\begin{flalign*}
\lambda_0=\frac{\Lambda_m\theta^2 }{\Lambda_h(\int_0^\infty \pi_h(s)\d s)^2} 
\left(\int_0^\infty \int_0^\infty  \beta_h(\tau+\xi,\tau) e^{-\int_0^\tau (\mu_h(\sigma+\xi)+\nu_h(\sigma+\xi,\sigma)+\gamma_h(\sigma+\xi,\sigma))\d \sigma}   \pi_h(\xi)  \d \xi~\d \tau\right)\\
\times\left(\int_0^\infty \int_0^\infty \beta_m(\tau+\xi,\tau)  e^{-\int_0^\tau (\mu_m(\sigma+\xi)+\nu_m(\sigma+\xi,\sigma))\d \sigma}  \pi_m(\xi) \d \xi~\d \tau\right).
\end{flalign*}
It also holds that
$$\tilde{H}_h(\mathbf{1}_{\Omega_h} \pi_h)=\lambda_0 \mathbf{1}_{\Omega_h}\pi_h$$
hence $\lambda_0$ is an eigenvalue of $\tilde{H}_h$ associated to the eigenfunction $\mathbf{1}_{\Omega_h}\pi_h$. It follows from a version of the Krein-Rutman theorem (see \textit{e.g.} \cite[Corollary 4.2.15, p. 273]{MeyerNieberg91}) that the spectral radius of $\tilde{H}_h$ is the only eigenvalue associated to a positive eigenfunction. Since $\pi_h>0$ a.e. on $\R_+$, we deduce that
$$r_\sigma(\tilde{H}_h)=\lambda_0.$$ 
Now, since $\lambda_0>0$, we see that $r_\sigma(H_h)>0$ by means of \eqref{Eq:spectral_rad>0}. It follows from \cite[Lemma 4.2.10, p. 269]{MeyerNieberg91} that there exists a positive eigenfunction $\phi_h\in L^1_+(\R_+)\setminus\{0\}$ such that 
$$H_h(\phi_h)=r_\sigma(H_h)\phi_h.$$
It follows that
$$\tilde{H}_h(\mathbf{1}_{\Omega_h}\phi_h)=r_\sigma(H_h)\mathbf{1}_{\Omega_h}\phi_h$$
so that $r_\sigma(H_h)>0$ is an eigenvalue of $\tilde{H}_h$ associated to $\mathbf{1}_{\Omega_h}\phi_h\in L^1_+(\Omega_h)\setminus \{0\}$. Again, by Krein-Rutman theorem, it follows that 
$$r_\sigma(H_h)=r_\sigma(\tilde{H}_h)=\lambda_0.$$
Using the same arguments, with $\pi_m$ instead of $\pi_h$, we can show that
$$r_\sigma(H_m)=\lambda_0.$$
Finally, using \eqref{Eq:R0_h}, it follows that
$$\RR_0=\sqrt{\lambda_0}$$
which ends the proof.

\subsection{Proof of Theorem \ref{Thm:stab}}

This section is dedicated to the stability analysis of the disease-free equilibrium of Model \eqref{Eq:Model}. 

Let $A_0$ be the part of $A$ in $X_0$, then denote by $\{T_{A_0}(t)\}_{t\geq 0}$ the positive semigroup generated by $A_0$. We know from Section \ref{Sec:Sg_formulation} that $(-\mu_0,\infty)\subset \rho(A_0)$ and consequently $s(A_0)\leq -\mu_0$ (where $s(A_0)$ denotes the spectral bound of $A_0$). Since the semigroup $\{T_{A_0}(t)\}_{t\geq 0}$ is positive, it follows that $\omega_0(\{T_{A_0}(t)\}_{t\geq 0})=s(A_0)$ (where $\omega_0$ denotes the growth bound) by using \cite[Theorem VI.1.15, p. 358]{EngelNagel2000}. Moreover, we know that $\omega_{\ess}(\{T_{A_0}(t)\}_{t\geq 0})$, the essential growth bound of $\{T_{A_0}(t)\}_{t\geq 0}$, satisfies $\omega_{\ess}(\{T_{A_0}(t)\}_{t\geq 0})\leq \omega_0(\{T_{A_0}(t)\}_{t\geq 0})$. We then have 
$$\omega_{\ess}(\{T_{A_0}(t)\}_{t\geq 0})\leq -\mu_0<0.$$

\textbf{Step 1:} we show that the operator $(D_{E_0}F_{\ep})_1$ is compact. We first rewrite it as:
$$(D_{E_0}F_{\ep})_1=(0,-G_1,G_1,0,0,0,0,-G_2,G_2,0)^T$$
where $G_1,G_2: X_{0}\to L^1_+(\R_+)$. It remains to show that both operators $G_1$ and $G_2$ are compact. To this end, we will use the classical Rietz-Fréchet-Kolmogorov (RFK) criterion in $L^1$ (see \textit{e.g} \cite[Theorem X.1, p. 275]{Yosida80}). Let $h\in\R_+$ and $Z\subset X_{0}$ be a bounded subset of $X_{0}$. Then there exists a positive constant $m>0$ such that $\|u\|_{X}\leq m$ for each $u\in Z$. Let $\TT_h$ be the translation operator in $L^1$, \textit{i.e.}
$$\TT_h(\phi)=\phi(\cdot+h)$$
for each $\phi\in L^1(\R_+)$. We have on one hand
$$
\|\TT_h(G_1(u))-G_1(u)\|_{L^1(\R_+)}\leq \left(\frac{m\theta \|\beta_m\|_{L^\infty} }{N^0_h}\right)\left\|\TT_h\left(S^0_h\right)-S^0_h\right\|_{L^1(\R_+)}\underset{h\to 0}{\to}0$$
since $S^0_h\in L^1(\R_+)$. It implies that
$$\sup_{u\in Z}\|\TT_h(G_1(u))-G_1(u)\|_{L^1(\R_+)}\underset{h\to 0}{\to}0.$$
On the other hand, by using Lebesgue theorem, we see that
$$\sup_{u\in Z}\int_c^\infty G_1(u)(a)\d a\leq \left(\frac{m\theta \|\beta_m\|_{L^\infty}}{N^0_h}\right)\int_c^\infty S^0_h(a)\d a\underset{c\to \infty}{\to}0.$$
The consequence of the two last points and the RFK criterion dwell in the compactness of $G_1(Z)$ in $L^1(\R_+)$. It proves that $G_1$ and $-G_1$ are compact. Similarly we can show that $G_2$ and $-G_2$ are also compact. It follows that the operator $(D_{E_0}F_{\ep})_1$ is compact.

\textbf{Step 2:} we now show that $\omega_{0}(\{T_{(A+(D_{E_0}F_{\ep})_2)_0}\}_{t\geq 0})\leq -\mu_0$, where $\{T_{(A+(D_{E_0}F_{\ep})_2)_0}(t)\}_{t\geq 0}$ is the $C_0$-semigroup generated by $(A+(D_{E_0}F_{\ep})_2)_0$, that is the part of $A+(D_{E_0}F_{\ep})_2$ in $X_0$. Let $\hat{u}_0\in X_{0}$. By means of the Volterra integral formulation stated in Theorem \ref{Thm:existence}, we can compute the expression of the latter semigroup, as follows:
$$T_{(A+(D_{E_0}F_{\ep})_2)_0}(t)\hat{u}_0=(0,S_h(t,\cdot),0,I_h(t,\cdot,\cdot),0,R_h(t,\cdot,\cdot),0,S_m(t,\cdot),0 ,I_m(t,\cdot,\cdot))^T$$
where
\begin{equation*}
\begin{array}{lll}
 S_h(t,a)=\left\{
\begin{array}{lll}
  \int_0^a \left(\int_0^\infty k_h(s,\eta)R_h(t+s-a,s,\eta)\d \eta\right)e^{-\int_{s}^a \mu_h(\xi)\d \xi}\d s, \quad \forall t>a, \\
\int_0^t \left(\int_0^\infty k_h(a+s-t,\eta) R_h(s,a+s-t,\eta)\d \eta \right)e^{-\int_s^t \mu_h(a+\xi-t)\d \xi}\d s \\
\qquad+S_{h,0}(a-t) e^{-\int_{a-t}^a \mu_h(s)\d s}, \quad \forall a\geq t,
    \end{array}
    \right.
\\
I_h(t,a,\tau)=\mathbf{1}_{\{a\geq\tau\geq t\}}\left(I_{h,0}(a-t,\tau-t)e^{-\int_{\tau-t}^\tau (\mu_h(s+a-\tau)+\nu_h(s+a-\tau,s)+\gamma_h(s+a-\tau,s))\d s}\right),
\\
R_h(t,a,\eta)=\left\{
\begin{array}{lll}
\left(\int_0^\infty \gamma_h(a-\eta,\tau)I_h(t-\eta,a-\eta,\tau)\d \tau\right)e^{-\int_0^\eta (\mu_h(s+a-\eta)+k_h(s+a-\eta,s))\d s}, \ \forall t>\eta, \ \forall a\geq \eta,  \vspace{0.1cm} \\
R_{h,0}(a-t,\eta-t)e^{-\int_{\eta-t}^\eta (\mu_h(s+a-\eta)+k_h(s+a-\eta,s))\d s}, \quad \forall a\geq \eta\geq t,
    \end{array}
    \right.
\\
S_m(t,a)=\mathbf{1}_{[t,\infty)}(a) S_{m,0}(a-t)e^{-\int_{a-t}^a \mu_m(s)\d s}, \quad \forall a\geq \eta\geq t,
\\
I_m(t,a,\tau)=\mathbf{1}_{\{a\geq \tau \geq t\}}I_{m,0}(a-t,\tau-t)e^{-\int_{\tau-t}^\tau (\mu_m(s+a-\tau)+\nu_m(s+a-\tau,s))\d s}.
    \end{array}
\end{equation*}
It follows that for each $t\geq 0$ we have the following inequalities
$$\|I_h(t,\cdot,\cdot)\|_{L^1(\R_+^2)}\leq \|I_{h,0}\|_{L^1(\R_+^2)}e^{-\mu_0 t}, \qquad \|S_m(t,\cdot)\|_{L^1(\R_+)}\leq \|S_{m,0}\|_{L^1(\R_+)}e^{-\mu_0 t}$$
and 
$$\|I_m(t,\cdot,\cdot)\|_{L^1(\R_+^2)}\leq \|I_{m,0}\|_{L^1(\R_+^2)}e^{-\mu_0 t}.$$
Moreover we see that
$$\int_t^\infty \int_{\eta}^\infty R_h(t,a,\eta)\d a~\d \eta\leq \|R_{h,0}\|_{L^1(\R_+^2)}e^{-\mu_0 t}$$
while
\begin{flalign*}
\int_0^t \int_{\eta}^\infty R_h(t,a,\eta)\d a~\d \eta&\leq \|\gamma_h\|_{L^\infty}\int_0^t e^{-\mu_0 \eta}\|I_h(t-\eta,\cdot,\cdot)\|_{L^1(\R_+^2)}\d \eta \\
&\leq \|\gamma_h\|_{L^\infty}\|I_{h,0}\|_{L^1(\R_+^2)}t e^{-\mu_0 t}
\end{flalign*}
It follows that
$$\|R_h(t,\cdot,\cdot)\|_{L^1(\R_+^2)}\leq \|R_{h,0}\|_{L^1(\R_+^2)}e^{-\mu_0 t}+\|\gamma_h\|_{L^\infty}\|I_{h,0}\|_{L^1(\R_+^2)}t e^{-\mu_0 t}$$
for each $t\geq 0$. Finally, we need the following estimates for $S_h$:
\begin{flalign*}
\int_t^\infty S_h(t,a)\d a&\leq \|S_{h,0}\|_{L^1(\R_+)}e^{-\mu_0 t}+\|k_h\|_{L^\infty}e^{-\mu_0 t}\int_{t}^\infty \int_0^t e^{\mu_0 s}\int_0^\infty R_h(s,a+s-t,\eta)\d \eta~\d s~\d a \\
&\leq  \|S_{h,0}\|_{L^1(\R_+)}e^{-\mu_0 t}+\|k_h\|_{L^\infty}e^{-\mu_0 t}\int_0^t \|R_h(s,\cdot,\cdot)\|_{L^1(\R_+^2)}e^{\mu_0 s}\d s \\
&\leq \left(\|S_{h,0}\|_{L^1(\R_+)}+t\|k_h\|_{L^\infty}\|R_{h,0}\|_{L^1(\R_+^2)}+\|k_h\|_{L^\infty}\|\gamma_h\|_{L^\infty}\|I_{h,0}\|_{L^1(\R_+^2)}\left(\frac{t^2}{2}\right)\right)e^{-\mu_0 t}
\end{flalign*}
and 
\begin{flalign*}
\int_0^t S_h(t,a)\d a&\leq \int_0^t \int_0^a e^{-\mu_0(a-s)}\int_0^\infty k_h(s,\eta)R_h(t+s-a,s,\eta)\d \eta~\d s~\d a \\
&\leq \|k_h\|_{L^\infty}\int_0^t \int_0^a e^{-\mu_0 \xi} \int_0^\infty R_h(t-\xi,a-\xi,\eta)\d \eta~\d \xi~\d a\\
&\leq \int_0^t \|R_h(t-\xi,\cdot,\cdot)\|_{L^1(\R_+^2)}e^{-\mu_0 \xi}\d \xi\\
&\leq \|R_{h,0}\|_{L^1(\R_+^2)}te^{-\mu_0 t}+\|\gamma_h\|_{L^\infty}\|I_{h,0}\|_{L^1(\R_+^2)}e^{-\mu_0 t}\left(\frac{t^2}{2}\right).
\end{flalign*}
We thus get
\begin{flalign*}
\|S_h(t,\cdot)\|_{L^1(\R_+)}\leq& \left(\|S_{h,0}\|_{L^1(\R_+)}+t\|k_h\|_{L^\infty}\|R_{h,0}\|_{L^1(\R_+^2)}+\|k_h\|_{L^\infty}\|\gamma_h\|_{L^\infty}\|I_{h,0}\|_{L^1(\R_+^2)}\left(\frac{t^2}{2}\right)\right)e^{-\mu_0 t} \\
&+\|R_{h,0}\|_{L^1(\R_+^2)}te^{-\mu_0 t}+\|\gamma_h\|_{L^\infty}\|I_{h,0}\|_{L^1(\R_+^2)}e^{-\mu_0 t}\left(\frac{t^2}{2}\right).
\end{flalign*}
It follows that there exists a positive constant $c>0$ such that
$$\left\|T_{(A+(D_{E_0}F_{\ep})_2)_0}(t)\right\|_{\L(X)}\leq c(1+t+t^2)e^{-\mu_0 t}$$
and
$$
\omega_0\left(\left\{T_{(A+(D_{E_0}F_{\ep})_2)_0}(t)\right\}_{t\geq 0}\right)=\lim_{t\to\infty}\frac{1}{t}\ln\left(\left\|T_{(A+(D_{E_0}F_{\ep})_2)_0}(t)\right\|_{\L(X)}\right)\leq -\mu_0.$$

\textbf{Step 3:} we now prove Theorem \ref{Thm:stab}. Since $(D_{E_0}F_{\ep})_1$ is a compact bounded operator by Step 1, then it follows that
$$\omega_{\ess}\left(\left\{T_{(A+D_{E_0}F_{\ep})_0}(t)\right\}_{t\geq 0}\right)=\omega_{\ess}\left(\left\{T_{(A+(D_{E_0}F_{\ep})_2)_0}(t)\right\}_{t\geq 0}\right)
$$
by using \cite[Theorem 1.2]{DucrotMagal2008}. It follows by Step 2 that
$$\omega_{\ess}\left(\left\{T_{(A+D_{E_0}F_{\ep})_0}(t)\right\}_{t\geq 0}\right)\leq \omega_{0}\left(\left\{T_{(A+(D_{E_0}F_{\ep})_2)_0}(t)\right\}_{t\geq 0}\right)\leq -\mu_0.$$
From \cite[Corollary IV. 2.11, p. 258]{EngelNagel2000} we deduce that
$$\{\lambda\in \sigma((A+D_{E_0}F_{\ep})_0), \Re(\lambda)\geq -\mu_0\}$$
is finite and composed (at most) of isolated eigenvalues with finite algebraic multiplicity, where $\sigma(\cdot)$ denotes the spectrum. Consequently, it remains to study the punctual spectrum of $(A+D_{E_0}F_{\ep})_0$. Using \cite[Proposition 4.19, p. 20]{Webb85}, we know that if $s((A+D_{E_0}F_{\ep})_0)<0$ then $E_0$ is locally asymptotically stable, while if $s((A+D_{E_0}F_{\ep})_0)>0$ then $E_0$ is unstable. We consider exponential solutions, \textit{i.e.} of the form $u(t)=e^{\lambda t}v$, with $0\neq v=(0,S_h,0,I_h,0,R_h,0,S_m,0,I_m)\in D(A)$ and $\lambda\in \C$. We obtain the following system:
\begin{equation*}
\left\{
\begin{array}{rcl}
S_h'(a)&=&\int_0^\infty k_h(a,\eta) R_h(a,\eta)\d \eta-\mu_h(a) S_h(a)-\frac{S^0_h(a)}{N^0_h}\int_0^\infty \int_0^\infty \theta \beta_m(s,\tau)I_m(s,\tau)\d s~\d \tau\\
\left(\frac{\partial}{\partial a}+\frac{\partial}{\partial \tau}\right)I_h(a,\tau)&=&-\left(\mu_h(a)+\nu_h(a,\tau)+\gamma_h(a,\tau)\right)I_h(a,\tau), \vspace{0.1cm} \\
\left(\frac{\partial}{\partial a}+
\frac{\partial}{\partial \eta}\right)R_h(a,\eta)&=&-(\mu_h(a)+k_h(a,\eta))R_h(a,\eta), \vspace{0.1cm} \\
S_m'(a)&=&-\mu_m(a)S_m(a)-\frac{S^0_m(a)}{N^0_h}\int_0^\infty \int_0^\infty \theta \beta_h(s,\tau)I_h(s,\tau)\d s~\d \tau, \\
\left(\frac{\partial }{\partial a}+\frac{\partial }{\partial \tau}\right)I_m(a,\tau)&=&-(\mu_m(a)+\nu_m(a,\tau))I_m(a,\tau), \\
\end{array}
\right.
\end{equation*}
with the boundary conditions:
\begin{equation*}
\left\{
\begin{array}{rclll}
S_h(0)&=&0, &S_m(0)=0,\\
I_h(a,0)&=&\frac{S^0_h(a)}{N^0_h}\int_0^\infty \int_0^\infty \theta \beta_m(s,\tau)I_m(s,\tau)\d s~\d \tau, &I_h(0,\tau)=0, \vspace{0.1cm} \\
R_h(a,0)&=&\int_0^\infty \gamma_h(a,\tau)I_h(a,\tau) \d \tau, &R_h(0,\eta)=0, \vspace{0.1cm} \\
I_m(a,0)&=&\frac{S^0_m(a)}{N^0_h}\int_0^\infty \int_0^\infty \theta \beta_h(s,\tau)I_h(s,\tau)\d s~\d \tau, &I_m(0,\tau)=0
\end{array}.
\right.
\end{equation*}
We then get
$$I_h(a,\tau)=S^0_h(a-\tau)\left(\frac{\int_0^\infty \int_0^\infty \theta \beta_m(s,\tau)I_m(s,\tau)\d s~\d \tau}{N^0_h}\right)e^{-\int_0^\tau (\lambda+\mu_h(s+a-\tau)+\nu_h(s+a-\tau,s)+\gamma_h(s+a-\tau,s))\d s}$$
and
$$I_m(a,\tau)=S^0_m(a-\tau)\left(\frac{\int_0^\infty \int_0^\infty \theta \beta_h(s,\tau)I_h(s,\tau)\d s~\d \tau}{N^0_h}\right)e^{-\int_0^\tau (\lambda+\mu_m(s+a-\tau)+\nu_m(s+a-\tau,s))\d s}.$$
It follows that
\begin{flalign*}
\int_0^\infty \int_0^\infty \beta_h(a,\tau)I_h(a,\tau)\d a~\d \tau=&\Lambda_h\int_0^\infty\int_0^\infty \beta_h(a,\tau)e^{-\int_0^{a}\mu_h(s)\d s}e^{-\int_0^\tau (\lambda+\nu_h(s+a-\tau,s)+\gamma_h(s+a-\tau,s))\d s}\d a~\d \tau \\
&\times \left(\frac{\int_0^\infty \int_0^\infty \theta \beta_m(a,\tau)I_m(a,\tau)\d s~\d \tau}{N_h^0}\right)
\end{flalign*}
and
\begin{flalign*}
\int_0^\infty \int_0^\infty \beta_m(a,\tau)I_m(a,\tau)\d a~\d \tau=&\Lambda_m \int_0^\infty \int_0^\infty \beta_m(a,\tau)e^{-\int_0^a \mu_m(s)\d s}e^{-\int_0^\tau (\lambda+\nu_m(s+a-\tau,s))\d s}\d a~\d \tau \\
&\times \left(\frac{\int_0^\infty \int_0^\infty \theta \beta_h(a,\tau)I_h(a,\tau)\d a~\d \tau}{N_h^0}\right)
\end{flalign*}
whence
\begin{flalign*}
1=&\left(\frac{\Lambda_m \Lambda_h \theta^2}{(N^0_h)^2}\right)\left(\int_0^\infty \int_0^\infty \beta_m(a,\tau)e^{-\int_0^a \mu_m(s)\d s}e^{-\int_0^\tau (\lambda+\nu_m(s+a-\tau,s))\d s}\d a~\d \tau\right)\\
&\times \left(\int_0^\infty \int_0^\infty \beta_h(a,\tau)e^{-\int_0^a \mu_h(s)\d s}e^{-\int_0^\tau (\lambda+\nu_h(s+a-\tau,s)+\gamma_h(s+a-\tau,s))\d s}\d a~\d \tau\right)
\end{flalign*}
which is equivalent to $g(\lambda)=1$ where $g:\C\to \R$ is defined by
\begin{flalign*}
g(\lambda)=&\left(\frac{\Lambda_m \theta^2}{\Lambda_h\left(\int_0^\infty e^{-\int_0^a \mu_h(s)\d s}\right)^2}\right)\left(\int_0^\infty \int_0^\infty \beta_m(a,\tau)e^{-\int_0^a \mu_m(s)\d s}e^{-\int_0^\tau (\lambda+\nu_m(s+a-\tau,s))\d s}\d a~\d \tau\right) \\
&\times \left(\int_0^\infty \int_0^\infty \beta_h(a,\tau)e^{-\int_0^a \mu_h(s)\d s}e^{-\int_0^\tau (\lambda+\nu_h(s+a-\tau,s)+\gamma_h(s+a-\tau,s))\d s}\d a~\d \tau\right).
\end{flalign*}
Since $g(0)=R_0:= \mathcal{R}_0^2$, it readily follows that if $R_0>1$ then there exists $\lambda>0$ such that $g(\lambda)=1$ whence $s((A+D_{E_0}F_{\ep})_0)>0$ and $E_0$ is unstable. On the contrary, if $R_0<1$ then for each $\lambda\in \C$ such that $\Re(\lambda)\geq 0$, we have $g(\lambda)<1$, whence $s((A+D_{E_0}F_{\ep})_0)<0$ and $E_0$ is locally asymptotically sable, which ends the proof of Theorem \ref{Thm:stab}.

\section{Existence and stability of an endemic equilibrium} \label{Sec:equilibria}

\subsection{Existence of an endemic equilibrium}

We know that any endemic equilibrium
$$E^*=\left(S^*_h,I^*_h, R^*_h, S^*_m, I^*_m\right)$$
must satisfy the following equations:
\begin{flalign*}
S^*_h(a)=&\int_0^a \left(\int_0^\infty k_h(s,\eta) R^*_h(s,\eta)\d \eta \right)e^{-\int_s^a \mu_h(\xi)\d \xi-(a-s)\left(\frac{\int_0^\infty \int_0^\infty 
\theta \beta_m(\xi,\tau)I^*_m(\xi,\tau)\d \xi~\d \tau}{N^*_h}\right)}\d s \\
\quad &+\Lambda_h e^{-\int_0^a \mu_h(s)\d s-a \left(\frac{\int_0^\infty \int_0^\infty 
\theta \beta_m(s,\tau)I^*_m(s,\tau)\d s~\d \tau}{N^*_h}\right)}, \\
S^*_m(a)=&\Lambda_m e^{-\int_0^a \left(\mu_m(s)+\frac{\int_0^\infty \int_0^\infty \theta \beta_h(\xi,\tau)I^*_h(\xi,\tau)\d \xi~\d \tau}{N^*_h}\right)\d s},
\end{flalign*}
\begin{equation}
\begin{array}{l}
I^*_h(a,\tau)=\left(\frac{S^*_h(a-\tau)\int_0^\infty\int_0^\infty \theta \beta_m(s,\xi)I^*_m(s,\xi)\d s~\d \xi}{N^*_h}\right)e^{-\int_0^\tau (\mu_h(s+a-\tau)+\nu_h(s+a-\tau,s)+\gamma_h(s+a-\tau,s))\d s}, \label{Eq:equil_Ih}
\end{array}
\end{equation}
\begin{equation}
\begin{array}{l}\label{Eq:Equil_Iv}
I^*_m(a,\tau)=\left(\frac{\Lambda_m \int_0^\infty \int_0^\infty \theta \beta_h(s,\xi)I^*_h(s,\xi)\d s~\d \xi}{N^*_h}\right)e^{\frac{(\tau-a)\int_0^\infty \int_0^\infty \theta \beta_h(s,\xi)I^*_h(s,\xi)\d s~\d \xi}{N^*_h}-\int_0^a \mu_m(s)\d s-\int_0^\tau \nu_m(s+a-\tau,s)\d s},
\end{array}
\end{equation}
for each $a\geq \tau$,
\begin{flalign}
R^*_h(a,\eta)=\left(\int_0^\infty \gamma_h(a-\eta,\tau)I^*_h(a-\eta,\tau)\d \tau\right)\exp\left(-\int_0^\eta (\mu_h(s+a-\eta)+k_h(s+a-\eta,s))\d s\right), \label{Eq:equil_rh}
\end{flalign}
for each $a\geq \eta$. We also have
$$N^*_h=\int_0^\infty S^*_h(a)\d a+\int_0^\infty \int_0^\infty I^*_h(a,\tau)\d a~\d \tau+\int_0^\infty \int_0^\infty R^*_h(a,\eta)\d a~\d \eta.$$
Let the change of variable:
\begin{equation}\label{Eq:change_mar}
i^*_h=I^*_h/N^*_h, \qquad s^*_h=S^*_h/N^*_h, \qquad r^*_h=R^*_h/N^*_h.
\end{equation}
Using \eqref{Eq:equil_Ih}-\eqref{Eq:change_mar}, we find that
\begin{flalign}
&\int_0^\infty \int_\tau^\infty \theta \beta_h(a,\tau)i^*_h(a,\tau)\d a~\d \tau=\left(\int_0^\infty \int_\tau^\infty \theta \beta_m(a,\tau)I^*_m(a,\tau)\d a~\d \tau\right) \label{Eq:equil_betahIh} \\
&\times \left(\int_0^\infty \int_\tau^\infty \theta  \beta_h(a,\tau)\left(\frac{s^*_h(a-\tau)}{N^*_h}\right)e^{-\int_0^\tau (\mu_h(s+a-\tau)+\nu_h(s+a-\tau,s)+\gamma_h(s+a-\tau,s))\d s}\d a~\d \tau\right). \nonumber
\end{flalign}
From \eqref{Eq:Equil_Iv}-\eqref{Eq:change_mar}, we deduce that
\begin{flalign}
&\int_0^\infty \int_\tau^\infty \beta_m(a,\tau)I^*_m(a,\tau)\d a~\d \tau=\Lambda_m \left(\int_0^\infty \int_\tau^\infty \beta_h(a,\tau)i^*_h(a,\tau)\d a~\d \tau \right)   \label{Eq:equil_betavIv} \\
&\quad\times\left(\int_0^\infty \int_\tau^\infty \theta \beta_m(a,\tau)e^{-\int_0^a \mu_m(s)\d s-\int_0^\tau \nu_m(s+a-\tau,s)\d s+(\tau-a)\int_0^\infty \int_\xi^\infty \theta \beta_h(s,\xi)i^*_h(s,\xi)\d s~\d \xi}\d a~\d \tau\right).\nonumber
\end{flalign}
A necessary condition for an endemic equilibrium to exist is
$$\int_0^\infty \int_\tau^\infty \beta_h(a,\tau)i^*_h(a,\tau)\d a~\d \tau\neq 0, \qquad \int_0^\infty \int_\tau^\infty \beta_m(a,\tau)I^*_m(a,\tau)\d a~\d \tau\neq 0.$$
Hence, by means of \eqref{Eq:change_mar}-\eqref{Eq:equil_betahIh}-\eqref{Eq:equil_betavIv}, we see that a necessary and sufficient condition is:
\begin{flalign}\label{Eq:equil_cond1}
1=&\left(\frac{\Lambda_m \theta^2}{N^*_h}\right)\left(\int_0^\infty \int_\tau^\infty  \beta_h(a,\tau)s^*_h(a-\tau)e^{-\int_0^\tau (\mu_h(s+a-\tau)+\nu_h(s+a-\tau,s)+\gamma_h(s+a-\tau,s))\d s}\d a~\d \tau\right) \\
&\times \left(\int_0^\infty \int_\tau^\infty \beta_m(a,\tau)e^{-\int_0^a \mu_m(s)\d s-\int_0^\tau \nu_m(s+a-\tau,s)\d s+(\tau-a)}\int_0^\infty \int_\xi^\infty \theta \beta_h(s,\xi)i^*_h(s,\xi)\d s~\d \xi\right)\d a~\d \tau. \nonumber
\end{flalign}
Now let $n_h^*(a)$ be the total population of humans with age $a$ at equilibrium, \textit{i.e.}
\begin{equation}\label{Eq:nh_def}
n_h^*(a)=S_h^*(a)+\int_0^a I_h^*(a,\tau)\d \tau+\int_0^a R_h^*(a,\eta)\d \eta.
\end{equation}
Then we see that:
$$(n^*_h)'(a)=-\mu_h(a) n_h^*(a)-\int_0^a \nu_h(a,\tau)I_h^*(a,\tau)\d \tau, \quad n_h(0)=\Lambda_h$$
whence
\begin{equation}\label{Eq:nh}
    n_h^*(a)=\Lambda_h e^{-\int_0^a \mu_h(s)\d s}-\int_0^a \left(\int_0^s \nu_h(s,\tau)I^*_h(s,\tau)\d \tau\right)e^{-\int_s^a \mu_h(\xi)\d \xi}\d s, \quad \forall a\geq 0.
\end{equation}
We deduce that
$$N_h^*=\int_0^\infty n^*_h(a)\d a=\Lambda_h \int_0^\infty e^{-\int_0^a \mu_h(s)\d s}\d a-\int_0^\infty \int_0^a \left(\int_0^s \nu_h(s,\tau)I^*_h(s,\tau)\d \tau\right)e^{-\int_s^a \mu_h(\xi)\d \xi}\d s~\d a.$$
Dividing the latter equation by $N^*_h$ and using \eqref{Eq:change_mar}, we get
\begin{equation}\label{Eq:Nh_equil}
    N^*_h=\frac{\Lambda_h \int_0^\infty e^{-\int_0^a \mu_h(s)\d s}\d a}{1+\int_0^\infty \int_0^a \left(\int_0^s \nu_h(s,\tau)i^*_h(s,\tau)\d \tau\right)e^{-\int_s^a \mu_h(\xi)\d \xi}\d s~\d a}.
\end{equation}
Moreover it follows from \eqref{Eq:nh_def} that
$$s^*_h(a)=\frac{n^*_h(a)}{N^*_h}-\int_0^a i^*_h(a,\tau)\d \tau-\int_0^a r^*_h(a,\eta)\d \eta$$
whence
\begin{flalign}
s^*_h(a)=&\left(\frac{1+\int_0^\infty \int_0^\zeta \left(\int_0^s \nu_h(s,\tau)i^*_h(s,\tau)\d \tau\right)e^{-\int_s^\zeta \mu_h(\xi)\d \xi}\d s~\d \zeta}{\int_0^\infty e^{-\int_0^\zeta \mu_h(s)\d s}\d \zeta}\right)e^{-\int_0^a \mu_h(s)\d s} \label{Eq:equil_sh} \\
&-\int_0^a i^*_h(a,\tau)\d \tau -\int_0^a \left(\int_0^\infty \gamma_h(\tau)i^*_h(a-\eta,\tau)\d \tau\right)e^{-\int_0^\eta (\mu_h(s+a-\eta)+k_h(s))\d s}\d \eta.\nonumber \\
&-\int_0^a\left(\int_0^s \nu_h(s,\tau)i^*_h(s,\tau)\d \tau\right)e^{-\int_s^a \mu_h(\xi) \d \xi}\d s \nonumber
\end{flalign}
By setting $R_0:= \mathcal{R}_0^2$, it follows that the necessary and sufficient condition \eqref{Eq:equil_cond1}, with $N^*_h, s^*_h$ respectively given by \eqref{Eq:Nh_equil}-\eqref{Eq:equil_sh}, becomes
\begin{equation}
\label{Eq:equil_cond2}
\begin{array}{ll}
1=&R_0\left(1+\int_0^\infty \int_0^a \left(\int_0^s \nu_h(s,\tau) i_h^*(s,\tau)\d \tau\right) \exp\left(-\int_s^a \mu_h(\xi)\d \xi\right)\d s~\d a\right)^2\\
&\times\left(\frac{\int_0^\infty \int_\tau^\infty \beta_m(a,\tau)\exp\left(-\int_0^a \mu_m(s)\d s-\int_0^\tau \nu_m(s+a-\tau,s)\d s+(\tau-a)\int_0^\infty \int_\xi^\infty \theta \beta_h(s,\xi)i^*_h(s,\xi)\d s~\d \xi\right)\d a~\d \tau}{\int_0^\infty \int_\tau^\infty \beta_m(a,\tau)\exp\left(-\int_0^a \mu_m(s)\d s-\int_0^\tau \nu_m(s+a-\tau,s)\d s\right)}\right)  \\
&-\left(\int_0^\infty \int_\tau^\infty \theta \beta_m(a,\tau)e^{-\int_0^a \mu_m(s)\d s-\int_0^\tau \nu_m(s+a-\tau,s)\d s+(\tau-a)\int_0^\infty \int_\xi^\infty \theta \beta_h(s,\xi)i^*_h(s,\xi)\d s~\d \xi}\d a~\d \tau\right) \\
&\times \Lambda_m\left(\frac{1+\int_0^\infty \int_0^a \left(\int_0^s \nu_h(s,\tau) i_h^*(s,\tau)\d \tau\right) \exp\left(-\int_s^a \mu_h(\xi)\d \xi\right)\d s~\d a}{\Lambda_h \int_0^\infty \exp\left(-\int_0^a \mu_h(s)\d s\right)\d a}\right) \\
&\times \int_0^\infty \int_{\tau}^\infty \theta \beta_h(a,\tau)\exp\left(-\int_0^\tau (\mu_h(s+a-\tau)+\nu_h(s+a-\tau,s)+\gamma_h(s+a-\tau,s))\d s\right)\\
&\times \left[\int_0^{a-\tau} i_h^*(a-\tau,s)\d s+\int_0^{a-\tau}\left(\int_0^s \nu_h(s,\tau)i^*_h(s,\tau)\d \tau\right)e^{-\int_s^{a-\tau} \mu_h(\xi) \d \xi}\d s\right. \\
&\left.+\int_0^{a-\tau} \left(\int_0^\infty \gamma_h(a-\tau-\eta,\tau)i^*_h(a-\tau-\eta,\tau)\d \tau\right)e^{-\int_0^\eta (\mu_h(s+a-\tau-\eta)+k_h(s+a-\tau-\eta,s))\d s} \right]\d a~\d \tau
\end{array}
\end{equation}
The latter condition has only $i_h^*$ as unknown, but is though difficult to exploit because $i_h^*$ depends both on $a$ and $\tau$.

\subsection{Proof of Theorem \ref{Thm:Bif}: bifurcations of an endemic equilibrium in a particular case}
Here, we assume that Assumptions \ref{Assump} and \ref{Assump-relax-human-age} hold. 
The threshold $\mathcal{R}_0^2$ rewrites as
\begin{flalign}\label{Eq:R0_2}
\mathcal{R}_0^2=&\left(\frac{\Lambda_m \mu_h \theta^2}{\Lambda_h}\right)\left(\int_0^\infty \int_\tau^\infty \beta_m(a,\tau)e^{-\int_0^a \mu_m(s)\d s-\int_0^\tau \nu_m(s+a-\tau,s)\d s}\d a~\d \tau\right)  \\
&\times \left(\int_0^\infty  \beta_h(\tau)e^{-\int_0^\tau (\mu_h+\nu_h(s)+\gamma_h(s))\d s}\d \tau\right). \nonumber
\end{flalign}
Let us recall the constant $C_{\bif}$ given by \eqref{Eq:constant_bif}. We can now give details on the proof of Theorem \ref{Thm:Bif} in several steps.
\textbf{Step 1:} we start by reminding that a theoretical necessary and sufficient condition to get an endemic equilibrium is given in the general case by \eqref{Eq:equil_cond1}. Using the change of variable \eqref{Eq:change_mar}, we see that the latter condition becomes:
\begin{flalign}
&1=\left(\frac{\Lambda_m s^*_h \theta^2}{N^*_h}\right)\left(\int_0^\infty\beta_h(\tau)e^{-\int_0^\tau (\mu_h+\nu_h(s)+\gamma_h(s))\d s}\d \tau\right) \label{Eq:equil_cond3}\\
&\left(\int_0^\infty \int_\tau^\infty \beta_m(a,\tau)e^{-\int_0^a \mu_m(s)\d s-\int_0^\tau \nu_m(s+a-\tau,s)\d s+(\tau-a)\int_0^\infty \theta \beta_h(s)i^*_h(s)\d s}\d a~\d \tau\right).\nonumber
\end{flalign}
From \eqref{Eq:change_mar} we get
\begin{equation}\label{Eq:somme=1}
s_h^*=1-\int_0^\infty i_h^*(\tau)\d \tau-\int_0^\infty r_h^*(\eta)\d \eta
\end{equation}
where 
$$\int_0^\infty r_h^*(\eta)\d \eta=\left(\int_0^\infty \gamma_h(s)i^*_h(s)\d s \right)\left(\int_0^\infty \exp\left(-\int_0^\eta (\mu_h+k_h(s))\d s\right)\d \eta\right)$$
is obtained from \eqref{Eq:equil_rh}. Moreover, using \eqref{Eq:Deriv_Nh}, we see that
$$\Lambda_h=\mu_h N^*_h+\int_0^\infty \nu_h(\tau)I_h^*(\tau)\d \tau$$
whence
$$N_h^*=\frac{\Lambda_h}{\mu_h+\int_0^\infty \nu_h(\tau)i^*_h(\tau)\d \tau}.$$
Using \eqref{Eq:R0_2} and setting $R_0:= \mathcal{R}_0^2$, it follows that the condition \eqref{Eq:equil_cond3} rewrites as
\begin{flalign}
1=&\left(\frac{R_0 \int_0^\infty \int_\tau^\infty \beta_m(a,\tau)e^{-\int_0^a \mu_m(s)\d s}e^{-\int_0^\tau \nu_m(s+a-\tau,s)\d s} e^{-(a-\tau)\int_0^\infty \theta \beta_h(s)i_h^*(s)\d s}\d a~\d \tau}{\int_0^\infty \int_\tau^\infty \beta_m(a,\tau)e^{-\int_0^a \mu_m(s)\d s} \d a~\d \tau}\right) \nonumber \\
\times&\left(1+\frac{\int_0^\infty \nu_h(\tau)i_h^*(\tau)\d \tau}{\mu_h}\right)\left[1-\left(\int_0^\infty i_h^*(\tau)\d \tau\right)-\left(\int_0^\infty \gamma_h(\tau)i_h^*(\tau)\d \tau\right)\int_0^\infty e^{-\int_0^\eta (\mu_h+k_h(s))\d s}\d \eta\right] \label{Eq:equil_cond4}
\end{flalign}
which may also be deduced from \eqref{Eq:equil_cond2}. We now write $i_h$ as:
\begin{equation}\label{Eq:i_h-K}
i_h^*(\tau)=K\exp\left(-\int_0^\tau (\mu_h+\nu_h(s)+\gamma_h(s))\d s\right)
\end{equation}
for each $\tau\in \R_+$, with 
\begin{equation}\label{Eq:K}
K=\frac{s_h^*\theta}{N_h^*}\int_0^\infty \int_\tau^\infty \beta_m(a,\tau) I_m^*(a,\tau)\d a~\d \tau.
\end{equation}
It is clear that, under Assumption \ref{Assump-relax-human-age}, an endemic equilibrium of \eqref{Eq:Model} exists if and only if $K>0$. Thus, the condition \eqref{Eq:equil_cond4} is equivalent to
$$f(R_0,K)=1$$
with $f$ defined by \eqref{Eq:function_f}.

It is then necessary that $K\in[0,\overline{K}]$ to have $f(R_0,K)\geq 0$ for any $R_0$, where $\overline{K}$ is given by: 
$$\overline{K}=\left[\int_0^\infty \exp\left(-\int_0^\tau (\mu_h+\nu_h(s)+\gamma_h(s))\d s\right)\left(1+\gamma_h(\tau)\int_0^\infty \exp\left(-\int_0^\eta (\mu_h+k_h(s))\d s\right)\d \eta \right)\d \tau\right]^{-1}.$$
We observe that for any $K\in[0,\overline{K}]$, it holds that $\int_0^\infty i_h^*(\tau)\d \tau\leq 1$, which is necessary by \eqref{Eq:somme=1}.

\textbf{Step 2:} we can now prove each item stated in the theorem.
\begin{enumerate}
\item Using the implicit function theorem, we get
$$\left({\frac{\d K}{\d R_0}}\right)_{\mid (K=0,R_0=1)}=-\frac{\partial_{R_0} f(0,1)}{\partial_K f(0,1)}.$$
On one hand we have
$$\partial_{R_0}f(0,1)=1$$
and on the other hand:
\begin{flalign}
&\partial_K f(R_0,K)=\left(\frac{R_0 \int_0^\infty \int_\tau^\infty \beta_m(a,\tau)e^{-\int_0^a \mu_m(s)\d s}e^{-\int_0^\tau \nu_m(s+a-\tau,s)\d s} e^{-(a-\tau)c_2 K}\d a~\d \tau}{\int_0^\infty \int_\tau^\infty \beta_m(a,\tau)e^{-\int_0^a \mu_m(s)\d s}e^{-\int_0^\tau \nu_m(s+a-\tau,s)\d s} \d a~\d \tau}\right) \label{Eq:derivf_K} \\
&\times\left[\frac{\int_0^\infty \nu_h(\tau)c_1(\tau)\d \tau}{\mu_h}\right.\left(1-2 K \int_0^\infty c_1(\tau)\d \tau -2K\left(\int_0^\infty \gamma_h(\tau)c_1(\tau)\d \tau\right)\int_0^\infty e^{-\int_0^\eta (\mu_h+k_h(s))\d s}\d \eta\right) \nonumber \\
&-\left.\left(\int_0^\infty c_1(\tau)\d \tau+\left(\int_0^\infty  \gamma_h(\tau)c_1(\tau)\d \tau\right)\int_0^\infty e^{-\int_0^\eta (\mu_h+k_h(s))\d s}\d \eta\right)\right]-\left(1+\frac{K\int_0^\infty \nu_h(\tau)c_1(\tau)\d \tau}{\mu_h}\right) \nonumber\\
&\times\left(\frac{c_2 R_0 \int_0^\infty \int_\tau^\infty \beta_m(a,\tau)(a-\tau)e^{-\int_0^a \mu_m(s)\d s}e^{-\int_0^\tau \nu_m(s+a-\tau,s)\d s}e^{-(a-\tau)c_2 K}\d a~\d \tau}{\int_0^\infty \int_\tau^\infty \beta_m(a,\tau)e^{-\int_0^a \mu_m(s)\d s}e^{-\int_0^\tau \nu_m(s+a-\tau,s)\d s} \d a~\d \tau}\right)  \nonumber \\
&\times\left(1-K\int_0^\infty c_1(\tau)\d \tau -K\left(\int_0^\infty \gamma_h(\tau)c_1(\tau)\d \tau\right)\int_0^\infty e^{-\int_0^\eta (\mu_h+k_h(s))\d s}\d \eta\right) \nonumber
\end{flalign}
with
$$c_1(\tau)= \exp\left(-\int_0^\tau (\mu_h+\nu_h(s)+\gamma_h(s))\d s \right), \qquad c_2=\theta \int_0^\infty \beta_h(\tau) c_1(\tau)\d \tau.$$
We then deduce that
\begin{flalign*}
\partial_K f(0,1)=&-\frac{c_2\int_0^\infty \int_\tau^\infty \beta_m(a,\tau) e^{-\int_0^a \mu_m(s)\d s}e^{-\int_0^\tau \nu_m(s+a-\tau,s)\d s}(a-\tau)\d a~\d \tau}{\int_0^\infty \int_\tau^\infty \beta_m(a,\tau) e^{-\int_0^a \mu_m(s)\d s}e^{-\int_0^\tau \nu_m(s+a-\tau,s)}\d a~\d \tau}\\
&+\int_0^\infty  c_1(\tau)\left(\frac{\nu_h(\tau)}{\mu_h}-1\right)\d \tau-\int_0^\infty \gamma_h(\tau)  e^{-\int_0^\tau(\mu_h+k_h(s))\d s}\d \tau.
\end{flalign*}
We know that a backward bifurcation occurs at $R_0=1$ if and only if $$\left({\frac{\d K}{\d R_0}}\right)_{\mid (K=0,R_0=1)}<0 $$
which is equivalent to $\partial_K f(0,1)>0$, \textit{i.e.} $C_\bif>0$.
\item Similarly, we know that a forward bifurcation occurs at $R_0=1$ if and only if
$$\left({\frac{\d K}{\d R_0}}\right)_{\mid (K=0,R_0=1)}>0$$
which amounts to $\partial_K f(0,1)<0$, \textit{i.e.} $C_\bif<0$.
\item For each $R_0$, we see that $f(R_0,0)=R_0$ and $\lim_{K\to \overline{K}}f(R_0,K)=0$. Consequently, if $R_0>1$, then there exists at least one solution $K$ to $f(R_0,K)=1$. Now, suppose that $R_0=1$ and $C_\bif>0$. Then on one hand we have $f(1,0)=1$ and on the other hand we have $\partial_K f(0,1)>0$ from the first point. Hence we have $f(1,\ep)>1$ for any $\ep>0$ small enough. Consequently there exists at least one $K>0$ such that $f(1,K)=1$.

\item Suppose that the condition \eqref{Eq:particular_case} holds. It follows from \eqref{Eq:derivf_K} that $\partial_K f(R_0,K)<0$ for each $R_0$, so that the function $K\longmapsto f(R_0,K)$ is strictly decreasing on $\R_+$. Since $f(R_0,0)=R_0$, we readily see that there exists an endemic equilibrium if and only if $R_0>1$, and in that case the equilibrium is unique. We can also note that a forward bifurcation occurs in this case since $C_\bif<0$ in this case.
\end{enumerate}

\section*{Funding}
QR and TL received support from the ANR STORM under grant agreement 16-CE35-0007.

\bibliography{Bibliography.bib}

\end{document}